\documentclass{article}
\usepackage[square,sort,comma,numbers]{natbib}
\usepackage{url}
\usepackage{authblk}

\usepackage[utf8]{inputenc}
\usepackage[T1]{fontenc}
\usepackage[english]{babel}

\usepackage{aligned-overset}

\usepackage{enumerate}
\usepackage{enumitem}
\usepackage{graphicx}
\usepackage{subcaption}
\usepackage{booktabs}
\usepackage{mathtools}
\usepackage{amsopn}
\usepackage{amsmath,bm}
\usepackage{amsthm}
\usepackage{amssymb}
\usepackage{amsthm}
\usepackage{amsfonts}
\usepackage{booktabs}
\usepackage{fancyhdr}
\usepackage{fourier}
\usepackage{nicefrac}
\usepackage{pgfplots}
\pgfplotsset{compat=1.18}
\usepackage{hyperref}
\usepackage{cleveref}

\usepackage[normalem]{ulem} 
\usepackage{makecell} 
\usepackage{arydshln} 

\usepackage{tikz}
\usepackage{tikz-cd}

\numberwithin{equation}{section}
\numberwithin{figure}{section}
\numberwithin{table}{section}

\usepackage[top=3cm, bottom=3cm,
  inner=3cm, outer=3cm]{geometry}			
\usepackage{float} 									
\usepackage[indent]{parskip}								
\usepackage{tcolorbox}

\tikzset{every picture/.style={line width=0.11mm}}
\newcommand{\oPerpSymbol}{\begin{tikzpicture}[scale=0.134]
    \draw (0,-0.5)--(0,1); \draw (-0.866,-0.5)--(0.866,-0.5);
    \draw (0,0) circle [radius=1];
\end{tikzpicture}}
\newcommand{\operp}{\mathbin{\raisebox{-1pt}{\oPerpSymbol}}}

\newcommand{\RR}{\mathbb{R}}

\newcommand{\NN}{\mathbb{N}}

\newcommand{\Tan}{\mathbb{T}}

\DeclareMathOperator{\op}{op}
\DeclareMathOperator{\Id}{Id}

\DeclareMathOperator{\Alt}{Alt}

\DeclareMathOperator{\Ima}{Im}
\DeclareMathOperator{\Ker}{Ker}
\DeclareMathOperator{\dive}{div}
\DeclareMathOperator{\crl}{curl}

\newcommand{\dirder}{\nabla_{\normal}^{(\ell)}}

\newcommand{\mesh}{\mathcal{T}_h}
\newcommand{\meshcut}{\mathcal{T}_h^{\text{cut}}}
\newcommand{\meshimm}{\mathcal{T}_h^{\text{im}}}

\newcommand{\mcT}{\mathcal{T}}

\newcommand{\mcE}{\mathcal{E}}

\newcommand{\mcF}{\mathcal{F}}
\newcommand{\mcL}{\mathcal{L}}

\newcommand{\mcA}{\mathcal{A}}
\newcommand{\mcJ}{\mathcal{J}}
\newcommand{\mcFstab}{\mcF^\pa_h}

\newcommand{\Omdh}{\Omega_{h}}

\newcommand{\Vh}[2]{V_{h}^{#1,#2}}
\newcommand{\hatVh}[2]{\hat{V}_{h}^{#1,#2}}
\newcommand{\Pih}[2]{\Pi_{h}^{#1,#2}}

\newcommand{\pa}{\partial}

\newcommand{\normal}{{\boldsymbol{n}}}

\theoremstyle{definition}
\newtheorem{assumption}{Assumption}

\newtheorem{remark}{Remark}

\theoremstyle{plain}
\newtheorem{thm}{Theorem}
\newtheorem{prop}{Proposition}
\newtheorem{cor}{Corollary}
\newtheorem{lemma}{Lemma}

\newcommand{\email}[1]{\href{mailto:#1}{#1}}

\begin{document}

\title{Ghost stabilisation for cut finite element exterior calculus}

\author[1]{Daniele A. Di Pietro}
\author[1,2]{J\'er\^ome Droniou}
\author[1]{Erik Nilsson}
\affil[1]{IMAG, CNRS, Montpellier, France, \email{daniele.di-pietro@umontpellier.fr}, \email{jerome.droniou@cnrs.fr}, \email{erik.nilsson@umontpellier.fr}}
\affil[2]{School of Mathematics, Monash University, Melbourne, Australia}

\maketitle

\begin{abstract}
We introduce the cut finite element method in the language of finite element exterior calculus, by formulating a stabilisation -- for any form degree -- that makes the method robust with respect to the position of the interface relative to the mesh. We prove that the $L^2$-norm on the physical domain augmented with this stabilisation is uniformly equivalent to the $L^2$-norm on the ``active'' mesh that contains all the degrees of freedom of the finite element space (including those external to the physical domain). We show how this CutFEEC method can be applied to discretize the Hodge Laplace equations on an unfitted mesh, in any dimension and any topology, and prove stability and optimal convergence. 
Two numerical examples are provided, with convergence and condition number scaling independent of the position of the boundary with respect to the background mesh. We first solve the Hodge Laplace equation on a conforming finite element space of $H^{\text{curl}}$ posed on a filled torus. The second example extends beyond the Hodge Laplace problem and illustrates the importance of stabilisation for pressure robustness in unfitted schemes, using an $H^{\text{curl}}$-formulation for the Stokes equations.
  \smallskip\\
  \textbf{MSC:} 65N30, 
  14F40 
  \smallskip\\
  \textbf{Key words:} Cut finite element method, FEEC, Hodge Laplace equation, Hilbert complex
\end{abstract}


\section{Introduction}

A successful strand of the literature on the numerical approximation of partial differential equations (PDEs) has emphasised the relevance of Hilbert complexes~\cite{bruning1992hilbert} in the design of stable numerical schemes.
In this context, differential forms provide a natural unifying language to treat mixed problems set in arbitrary space dimension, possibly on manifolds, and where a combination of differential operators appear~\cite{Bossavit:02,Hiptmair:02,arnold_feec_homological_2006,arnold_feec_hodge_2010,diPDr2023exteriorpolyt,Di-Pietro.Droniou.ea:25*1}.
Finite Element Exterior Calculus (FEEC), particularly relevant for the present work, provides a comprehensive view of full and trimmed finite element approximations of the de Rham complex; see the monograph~\cite{arnold_feec_2018} for an introduction.
At the same time, the Cut Finite Element Method (CutFEM) has emerged as a versatile technique for handling PDEs on domains with complex geometries,
removing the need for the mesh to be compliant with the domain's boundary; see the recent review \cite{burman2025cutfemreview} and references therein.
Through weak imposition of boundary conditions and carefully constructed stabilisation terms over facets close to the boundary, CutFEM leads to stable and well-conditioned discretisations.
Moreover, such properties hold independently of the boundary position relative to the unfitted mesh.
This work aims at merging these research avenues by proposing a Cut Finite Element Exterior Calculus (CutFEEC) framework that adapts the CutFEM techniques to FEEC.
This merger has notably been made possible by the recent extension of the CutFEM technology to mixed problems~\cite{frachon2024darcy,frachon2024stokes,frachon2024lagrange}, which are a natural application of FEEC~\cite{arnold_feec_2018}.

The prime example of a Hilbert complex is the de Rham complex, which, for an open connected domain $\Omega \subset \mathbb{R}^3$, reads
\[
0\hookrightarrow  H^1(\Omega) \xrightarrow{\nabla} \boldsymbol{H}^{\crl}(\Omega) \xrightarrow{\crl} \boldsymbol{H}^{\dive}(\Omega) \xrightarrow{\dive} L^2(\Omega) \to 0.
\]
Using the language of differential forms, the de Rham complex can be extended to a domain $\Omega$ of any dimension $n\in\NN$ as:
\begin{equation}\label{eq:continuous.complex}
0\hookrightarrow  H\Lambda^0\Omega \xrightarrow{d^0} H\Lambda^1\Omega \xrightarrow{d^1} \dots \xrightarrow{d^{n-1}} H\Lambda^n\Omega \to 0,
\end{equation}
where $d^k$ is the exterior derivative and $H\Lambda^k\Omega$ is the space of $L^2$-integrable differential $k$-forms on $\Omega$ with $L^2$-integrable exterior derivative.
The de Rham complex enters the well-posedness analysis of PDEs through its cohomology spaces
\[
\mathfrak{H}^k\coloneq \Ker d^k / \Ima d^{k-1}.
\]
These spaces relate to the topology of the domain via their dimension. Preserving these homological structures at the discrete level is crucial for the stability of numerical schemes. 
A paradigmatic example of a PDE problem whose well-posedness hinges on the de Rham complex is the Hodge Laplace equation; see, e.g.,~\cite[Chapter~4]{arnold_feec_2018}.
Its mixed formulation naturally leads to a mixed problem where the unknowns are the codifferential, the (scalar- or vector-) potential, and a Lagrange multiplier to enforce $L^2$-orthogonality with respect to harmonic forms.

The Finite Element (FE) approximation of mixed problems became a research topic starting from the late 1970s, when the first conforming approximations of vector-valued spaces in the de Rham complex were developed~\cite{nedelec1980mixed,raviart1977mixed,bossavit1988rationale}. 
Traditional FE methods, such as the ones cited above, require meshes that conform to the geometry of the domain. These can be challenging to generate for intricate geometries or in the presence of evolving interfaces. CutFEM addresses this problem by permitting the interface to cut through the elements of a background mesh, thereby simplifying mesh generation and adaptation \cite{burman2015cutfem}. However, ensuring stability and accuracy in such unfitted methods necessitates appropriate stabilisation techniques \cite{Bu10Ghost}. Otherwise, when the geometry cuts the mesh in especially nasty ways ($|T\cap\Omega|\ll |T|$ for an element $T$ in the mesh), one observes a significant degradation of the condition number of the associated linear system. 
Recent works have addressed the application of CutFEM techniques to mixed problems.
In \cite{frachon2024darcy,frachon2024lagrange}, the authors introduced a stable CutFEM for Darcy flow problems which, thanks to conformity with respect to the tail end of the de Rham complex, preserves the divergence of the velocity field.
Using a similar approach, one can also preserve the divergence condition of Stokes flow in an unfitted setting~\cite{frachon2024stokes}.

In this work, we develop a new CutFEM stabilisation to design unfitted robust discrete $L^2$-products for any space in the de Rham complex of differential forms.
The stabilisation can be seen as a generalisation of the mixed ghost penalty term introduced in \cite{frachon2024darcy,frachon2024stokes,nilsson2024generalized}.
Orthogonalising against the harmonic forms using an $L^2$-product including this stabilisation, we construct an arbitrary-order numerical method which is unfitted but still compatible with the continuous de Rham complex, thus ensuring its natural stability and preservation of relevant constraints (e.g., zero divergence of electric field in Maxwell's equations). 
We develop facet-based stabilisation terms, but we mention that there exist other equivalent stabilisation terms in the CutFEM literature, though these have not been developed for differential forms, see \cite[Section 4.2]{burman2025cutfemreview}. 
We note that the commonly used assumption of a quasi-uniform mesh \cite{burman2025cutfemreview} is not needed in our construction or analysis. It is sufficient for the relevant arguments to use a local quasi-uniformity, which does not entail particular restrictions on the mesh. 
In \cite{massing2014stabilized}, the quasi-uniform assumption is used but a similar remark is made \cite[Remark 2.3]{gurkan2019stabilized}.

A particularly relevant novelty results from applying this general construction to the special case of $1$-forms, which is key to
discretize problems such as Maxwell's equations or the $H^{\crl}$-elliptic problem.
CutFEM has, indeed, seen little development for problems involving $H^{\crl}$; see the recent and exhaustive review \cite{burman2025cutfemreview}.
There are some previous works which use an unfitted discontinuous Galerkin approach \cite{chen2023highorderunfittedmaxwell-DG,li2023reconstructedHcurl-elliptic-DG}.
Otherwise, works which do not use discontinuous/broken finite element spaces are essentially limited to
\cite{yang2024unfittedmaxwell-Pk,wang2025mixedHcurl-elliptic-Ned,frachon2024stokes}\cite[Paper D]{nilsson2024generalized}.
Of these, \cite{frachon2024stokes} and \cite[Paper D]{nilsson2024generalized} do not contain a stability analysis; \cite{frachon2024stokes} considers a vorticity-velocity-pressure formulation for Stokes flow wherein the vorticity is discretised by a $H^{\crl}$-conforming variable, while \cite[Paper D]{nilsson2024generalized} considers numerical experiments of a set of standard formulations of Maxwell's equations made unfitted.
The works \cite{yang2024unfittedmaxwell-Pk,wang2025mixedHcurl-elliptic-Ned} focus, respectively, on the time-harmonic Maxwell's equations and the $H^{\crl}$-elliptic interface problem.
In \cite[Remark 3]{yang2024unfittedmaxwell-Pk}, the authors mention that their proposed method works numerically when using nodal Lagrange elements, but error analysis is provided only for the case of using piecewise discontinuous polynomials.
The numerical scheme in \cite{ang2025mixedHcurl-elliptic-Ned} uses the lowest order $H^{\crl}$-conforming Nédélec elements of the second kind, but adds diagonal stabilisation terms in the Lagrange multiplier block \cite[(4.8)--(4.11)]{wang2025mixedHcurl-elliptic-Ned}, thus perturbing the saddle point nature of the problem. Similar stabilisation terms also appear in \cite{yang2024unfittedmaxwell-Pk}.
Moreover, in contrast to our approach, no work mentioned above considers the presence of non-trivial harmonic forms, which appear when considering domains of general topology.

The remainder of this paper is organised as follows.
In Section \ref{sec:continuous_hodge_laplace} we discuss a motivating example of a mixed problem, the Hodge Laplace problem, whose well-posedness relies on the properties of the de Rham complex, and we summarise some notions of FEEC.
In Section \ref{sec:unf_mesh} we recall the notion of unfitted mesh, and we then introduce in Section \ref{appendix:bdry_forms} the concepts of tangential and normal parts of differential forms, relative to a facet of the mesh. In Section \ref{sec:ghost_penalty} we introduce the ghost penalty stabilisation operator, and show the uniform equivalence between the $L^2$-norm on the physical domain and the $L^2$-norm on the extended domain containing all the degrees of freedom.
In Section \ref{sec:ingredients} we collect the remaining tools needing for proper error analysis.
In Section \ref{sec:CutFEEC} we provide the numerical scheme and prove stability, consistency, and convergence. 
Finally, in Section \ref{sec:numerics} we validate our results on a numerical example on the filled torus, demonstrating that the use of the stabilisation drastically reduces the condition number of the discrete system. Deviating from the Hodge Laplace equation with a second example, we also show that stabilisation from the point of view of equivalent norms is helpful in devising unfitted schemes for the Stokes that are pressure robust.


\section{A motivating example}\label{sec:continuous_hodge_laplace}

In this section we briefly discuss the Hodge Laplacian, an archetypal problem whose well-posedness hinges on the properties of the de Rham complex.
This example serves as a motivation for the development of a CutFEEC theory, and will be used in Section~\ref{sec:numerics} to numerically demonstrate the efficiency of the ghost penalty stabilisation of Section~\ref{sec:ghost_penalty}.

\subsection{Hodge Laplace equation}\label{sec:hodge.laplace}

Let $B\subset\RR^n$ be an open set.
We denote by $\Lambda^k B$ the space of smooth differential $k$-forms on $B$.
We introduce the $L^2$-product $(\omega,\zeta)_B \coloneqq \int_B \omega\wedge\star \zeta$ on $\Lambda^k B$,
where $\wedge$ denotes the wedge product and $\star$ the Hodge star operator \cite[(6.2)]{arnold_feec_2018}.
The $L^2$-norm over $B$ is denoted by $\|{\bullet}\|_B\coloneqq \sqrt{(\cdot,\cdot)_B}$. We denote by $L^2\Lambda^k B$ the space of $k$-forms over $B$ that are bounded with respect to this norm.
With the exterior derivative $d^k: L^2\Lambda^kB \to L^2\Lambda^{k+1}B$, viewed as an unbounded operator, we also define
\begin{align*}
  H\Lambda^k B &\coloneqq  \left\{ \omega\in L^2\Lambda^k B: d^k\omega\in L^2\Lambda^{k+1} B\right\},
\end{align*}
with norm
\begin{align*}
  \|\omega\|^2_{HB} \coloneqq  \|\omega\|^2_B+\|d^k\omega\|^2_B.
\end{align*}

Let now $\Omega\subset\RR^n$ be a connected and open subset of $\mathbb{R}^n$ with Lipschitz boundary $\partial \Omega$ and outward unit normal vector field $\normal$. 
Given $f\in L^2\Lambda^k\Omega$, the Hodge Laplace equation consists in seeking $\eta \in \Lambda^k \Omega$ such that
\begin{equation}\label{eqs:HodgeLaplace}
  \begin{aligned}
    (d^{k-1}\delta^k + \delta^{k+1} d^k)\eta &= f-\pi f &\qquad& \text{in $\Omega$,}
    \\
    \pi \eta &= 0 &\qquad& \text{in $\Omega$},
  \end{aligned}
\end{equation}
where the codifferential $\delta^{k+1}: L^2\Lambda^{k+1}\Omega \to L^2\Lambda^k\Omega$ is the adjoint of $d^k$ and
 $\pi\coloneqq \pi_{\mathfrak{H}^k}$ is the orthogonal projection onto the space of harmonic forms
\begin{align*}
  \mathfrak{H}^k \coloneqq \left\{
  \rho \in L^2 \Lambda^k\Omega: d^k \rho = 0 \text{ and } \delta^k \rho=0
  \right\} = \Ker d^k \cap \Ker \delta^k \cong \Ker d^k / \Ima d^{k-1}.
\end{align*}
Here $\cong$ means vector space isomorphism.
The boundary conditions are implicit in the definition of the Hodge Laplace operator $d^{k-1}\delta^k + \delta^{k+1} d^k$, and are revealed upon choosing the domain of the exterior derivative. For ease of presentation, let us in this work choose the domain
\begin{align*}
  D(d^k) &= H\Lambda^k\Omega,
\end{align*}
so that the domain of the codifferential becomes (by duality) the space of $L^2(\Omega)$-integrable forms $\omega$ with $\delta \omega \in L^2(\Omega)$ and $\gamma(\star \omega) =0.$
After recalling vector proxies in three space dimensions in Table~\ref{tab:vector_proxies3d}, in
Table~\ref{tab:hodgelaplace_ex}, we express problem~\eqref{eqs:HodgeLaplace} in terms of vector proxies for different values of $k$; cf.~\cite[Section~4.2]{arnold_feec_hodge_2010} for further details.

\begin{table}
  \centering
  \begin{tabular}{cccc}
    \toprule
    $k$
    & Form
    & Proxy
    & Sobolev space
    \\ \midrule
    0
    & $q$
    & $q$
    & $H^1(\Omega)$
    \\
    1
    & $v_1 dx^1 + v_2 dx^2 + v^3 dx^3$
    & $\boldsymbol{v} = (v_1,v_2,v_3)$
    & $\boldsymbol{H}^{\crl}(\Omega)$
    \\
    2
    & $w_1 dx^2\wedge dx^3 - w_2 dx^1\wedge dx^3 + w_3 dx^1\wedge dx^2$
    & $\boldsymbol{w} = (w_1,w_2,w_3)$
    & $\boldsymbol{H}^{\dive}(\Omega)$
    \\
    3
    & $r dx^1\wedge dx^2\wedge dx^3$
    & $r$
    & $L^2(\Omega)$
    \\
    \bottomrule
  \end{tabular}
  \caption{Vector proxies in $\RR^3$.}
  \label{tab:vector_proxies3d}
\end{table}

\begin{table}
  \centering
  \begin{tabular}{cccc}
    \toprule
    $k=0$ & $k=1$ & $k=2$ & $k=3$ \\
    \midrule
    \makecell[l]{
      $-\Delta q = f-\bar f,$\\
      $\nabla q\cdot \normal = 0,$\\
      $\bar q = 0$
    } &
    \makecell[l]{
      $\crl\crl \boldsymbol{v} - \nabla\dive \boldsymbol{v} = \boldsymbol{f} - \boldsymbol{\pi} \boldsymbol{f},$ \\
      \text{$\boldsymbol{v}\cdot \normal = 0$ and $\crl \boldsymbol{v} \times \normal = 0,$} \\
      $\boldsymbol{v} \perp \mathfrak{H}^1$
    } &
    \makecell[l]{
      $\crl\crl \boldsymbol{w} - \nabla\dive \boldsymbol{w} = \boldsymbol{f} - \boldsymbol{\pi} \boldsymbol{f},$ \\
      \text{$\boldsymbol{w} \times \normal = 0$ and $\dive \boldsymbol{w} = 0,$} \\
      $\boldsymbol{w} \perp \mathfrak{H}^2$
    } &
    \makecell[l]{
      $-\Delta r = f,$\\
      $r = 0$
    } \\
    \bottomrule
  \end{tabular}
  \caption{The Hodge Laplace equation in $\RR^3$ for different differential form orders. Entries are labeled $A$, $B$, $C$ where $A$ are the bulk equation(s) in $\Omega$, $B$ the boundary condition(s) on $\partial \Omega$, and $C$ is the cohomological compatibility condition. }
  \label{tab:hodgelaplace_ex}
\end{table}

For convenience of notation, we shall for the remainder omit the index and just write $d$ for the exterior derivative.
A variational mixed formulation of \eqref{eqs:HodgeLaplace} is (cf. \cite[Theorem~4.7]{arnold_feec_2018}):
Find $(\sigma,\eta,\lambda)\in H\Lambda^{k-1}\Omega \times H\Lambda^k\Omega \times \mathfrak{H}^k$ such that
\begin{equation}\label{eqs:mixed}
  \begin{aligned}
    (\sigma, \tau)_\Omega - (\eta, d\tau)_\Omega &= 0 \qquad &&\forall \tau\in H\Lambda^{k-1}\Omega, \\
    (d \sigma, \zeta)_\Omega + (d\eta,d\zeta)_{\Omega} + (\lambda, \zeta)_\Omega &= (f, \zeta)_\Omega \qquad &&\forall \zeta\in H\Lambda^{k}\Omega,  \\
    (\eta, \rho)_\Omega &= 0 \qquad &&\forall \rho\in \mathfrak{H}^{k}.
  \end{aligned}
\end{equation}
The above formulation enforces in a weak manner $\sigma = \delta \eta$ and $\lambda = \pi f$. The (natural) boundary conditions are $\gamma \star \eta = 0$ and  $\gamma \star d\eta = 0$.

\begin{remark}[Boundary conditions]
We emphasise that natural boundary conditions are particularly ame\-nable for unfitted methods such as CutFEM, since such conditions are handled as boundary integrals in the right hand side. In contrast, essential boundary conditions cannot be imposed strongly in the test and trial spaces like what is usually done in the fitted setting, and one must instead enforce them weakly via for example Nitsche's method \cite{burman2012cutfemnitsche}. This implies a discrepancy in the boundary de Rham complexes between the continuous and discrete formulations, which is an additional error to consider. As such, we leave this for future work in order to facilitate the clarity of the present work.
\end{remark}

Using crucial properties of the de Rham complex such as Hodge decompositions and Poincaré inequalities, it can be proved that this mixed formulation is well-posed \cite[Theorems~4.8 and~4.9]{arnold_feec_2018}.
Preserving these properties at the discrete level, as is the case for the compatible discretisations considered here, naturally leads to numerical schemes for which stability can be proved mimicking the arguments for well-posedness of the continuous problem.

\subsection{Principles of FEEC}

Assume here that $\Omega$ is a polytopal domain, and let $\mcT_h$ be a conforming simplicial mesh of $\Omega$.
The Finite Element Exterior Calculus (FEEC) method provides the tools for a discrete formulation of \eqref{eqs:mixed} which is well-posed \cite[(5.3), Section 5.2.3]{arnold_feec_2018}, but non-conforming in general due to the non-inclusion of the discrete space of harmonic forms into the continuous one (these spaces have the same dimension, so the former being included in the latter would actually mean that one could compute harmonic forms exactly).
We make a brief summary of the results which pertain to this work.

Let $r\ge 1$ be a polynomial degree. For each $k\in \{0,\dots,n\}$, let $\Vh{k}{r}\subset H\Lambda^k\Omega$ be a finite-dimensional space of $k$-forms that are piecewise polynomial of degree at most $r$ on $\mcT_h$ and such that $d\Vh{k}{r}\subset \Vh{k+1}{r-1}$.
We adopt the convention that $V_h^{-1,r}\coloneqq \{0\}$ and $V_h^{n+1,r}\coloneqq \{0\}$.
These spaces define a sub-complex of the de Rham complex:
\begin{equation}\label{eq:discrete.complex}
  \begin{tikzcd}
    0 \arrow[r,hook] & \Vh{0}{r} \arrow{r}{d} & \Vh{1}{r-1} \arrow{r}{d} & \cdots \arrow{r}{d} & \Vh{n}{r-n} \arrow{r} & 0.
  \end{tikzcd}
\end{equation}
It can be shown \cite{arnold_feec_2018}, using, e.g., bounded cochain maps between the continuous and discrete de Rham complexes, that the cohomology spaces of \eqref{eq:discrete.complex} are isomorphic to those of \eqref{eq:continuous.complex}.
This diagram (in which the polynomial degree decreases along the sequence) corresponds to the case of full finite element spaces, $\Vh{k}{r}\coloneqq P_r\Lambda^k\Omdh$. When using trimmed spaces, $P^-_r\Lambda^k\Omdh\subset P_r\Lambda^k\Omdh$ in the notation of \cite{arnold_feec_2018}, the degree $r$ remains the same along the diagram. For application of our results to the trimmed spaces, simply replace any occurrence of $\Vh{k}{\bullet}$ with $P^-_r\Lambda^k\Omdh$.

The discrete exterior derivative acting on a discrete space is just the restriction of the exterior derivative from the corresponding continuous space. When relevant, we shall denote it by $d_h\coloneqq  d|_{\Vh{k}{r}}$, otherwise by $d$ when there is no confusion.
Then, one defines the space of discrete harmonic forms as
\begin{equation}\label{eq:Hhk}
  \mathfrak{H}_h^k \coloneqq  \Ker d_h / d\Vh{k-1}{r+1} \cong \left\{\rho_h \in \Vh{k}{r} :\
  \text{$d\rho_h=0$ and $(\rho_h,d\tau_h)_\Omega = 0$ for all $\tau_h\in \Vh{k-1}{r+1}$}\right\}.
\end{equation}
The non-inclusion $\mathfrak{H}_h^k \not\subset \mathfrak{H}^k$ is due to forms in $\mathfrak{H}^k$ needing to be orthogonal to all of $dH\Lambda^{k-1}\Omega$, while being in $\mathfrak{H}_h^k$ only guarantees orthogonality to the subspace $d\Vh{k-1}{r+1}$.

As mentioned in Section \ref{sec:hodge.laplace}, the Hodge decomposition of the complex spaces is an essential tool to analyse both the continuous Hodge Laplace equation and its numerical discretisations. For the complex \eqref{eq:discrete.complex}, we have the following discrete Hodge decomposition:
\begin{align}\label{eq:discrete_hodge_decomp}
  \Vh{k}{r} = (\Ker d_h)^\perp \operp d\Vh{k-1}{r+1} \operp \mathfrak{H}_h^k.
\end{align}

Given $f\in L^2\Lambda^k\Omega$ the FEEC scheme for problem~\eqref{eqs:mixed} reads:
Find $(\sigma_h,\eta_h,\lambda_h)\in \Vh{k-1}{r+1} \times \Vh{k}{r} \times \mathfrak{H}_h^k$ such that
\begin{alignat*}{2}
    (\sigma_h, \tau_h)_\Omega - (\eta_h, d\tau_h)_\Omega &= 0  &&\qquad\forall \tau_h\in \Vh{k-1}{r+1}, \\
     (d \sigma_h, \zeta_h)_\Omega + (d\eta_h,d\zeta_h)_{\Omega} + (\lambda_h, \zeta_h)_\Omega &= (f, \zeta_h)_\Omega  &&\qquad\forall \zeta_h\in \Vh{k}{r}, \\
    (\eta_h, \rho_h)_\Omega &= 0  &&\qquad\forall \rho_h\in \mathfrak{H}_h^{k}.
\end{alignat*}
The scheme is well-posed and optimally convergent in the $H\Lambda^k\Omega$-norm \cite[Theorem 5.5]{arnold_feec_2018}.


\section{Unfitted mesh}\label{sec:unf_mesh}

Let $\mesh = \{ T\}$ be an \emph{active mesh} such that each (open) simplex $T \in \mesh$ has a nonempty intersection with $\Omega$, and for which
\begin{equation*}
  \Omega \subset
  \Omdh \coloneq \left(\bigcup_{T\in\mesh} \overline{T}\right)^\circ.
\end{equation*}
The open set $\Omega_h$ is hereafter referred to as the active domain.

\begin{figure}
  \centering
  \includegraphics[width=0.3\textwidth]{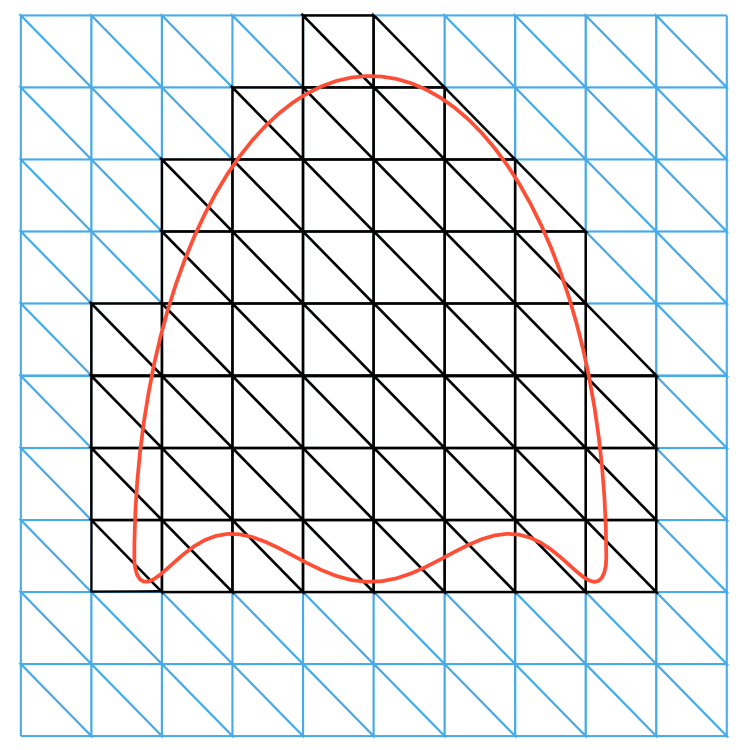}
  \caption{Illustration of the active mesh $\mcT_h$. The boundary $\partial\Omega$ of the domain is in red. The ``remainder'' of background mesh $\mcT_{0,h}\setminus \mesh$ is shown in light blue, while the active mesh $\mesh$ is shown in black.}
  \label{fig:active_mesh}
\end{figure}

\begin{remark}[Construction of active mesh]
  The term active mesh refers to the following construction.
  Define a background domain $\Omega_0\supset \Omega$, which is polytopal of dimension $n$.
  Consider a conforming simplicial mesh $\mcT_{0,h}$ of $\Omega_0$. Then, the active mesh is defined as
  $\mcT_h \coloneqq  \left\{ T\in\mcT_{0,h} : T\cap \Omega \neq \emptyset \right\}.$
  This construction is illustrated in Figure~\ref{fig:active_mesh}.
\end{remark}

We say $T\in\mesh$ is \emph{cut} if $T \not\subset \Omega$ and \emph{fully immersed} in $\Omega$ otherwise. We denote the set of cut elements by
\begin{equation*}
  \meshcut \coloneqq  \left\{
  T\in\mcT_h : T \not\subset \Omega
  \right\}.
\end{equation*}
We denote by $\mcFstab$ the set of stabilisation facets \cite[Figure 1]{frachon2024stokes}, which are facets internal to $\Omega_h$ for which at least one neighbouring element is intersected by $\partial\Omega$:
\[
\mcFstab \coloneqq \left\{
F \in \bigcup_{T\in\meshcut}\mcF_T \;:\; F \not\subset \partial \Omega_h
\right\},
\]
where $\mcF_T$ gathers the ($(n-1)$-dimensional) facets of $T \in \mesh$.
To each facet we associate a unique normal vector $\normal$.
The context will make it clear which normal is meant, that of a facet or that of $\partial\Omega$.
\begin{assumption}[Mesh regularity]\label{assu:mesh}
  For any $T \in \mcT_h$, let $\varrho_T$ and $h_T$ be, respectively, the diameter of the largest ball contained in $T$ and that of $T$.
  We assume the following:
  \begin{enumerate}[label=(\roman*)]
  \item (Shape regularity) For some constant $\kappa>0$ independent of the mesh size $h \coloneqq \max_{T \in \mesh} h_T$, we have
    \begin{equation}\label{eq:shape.regularity}
      \frac{h_T}{\varrho_T} \leq \kappa \qquad \forall T\in\mesh.
    \end{equation}
  \item\label{item:cut-to-uncut} (Uniformly bounded cut-to-uncut path) There exists an integer $N > 0$ independent of $h$ such that, for all $T \in \meshcut$, there exists a sequence of $N_T \le N$ elements $\{T_1 = T, T_2,\ldots,T_{N_T}\}$ such that $T_i$ and $T_{i+1}$ share a facet for all $1\le i\le N_T - 1$, and $T_{N_T}\in\mesh\setminus \meshcut$ is fully immersed in $\Omega$.
  \item For all $h<1$, the domain $\Omdh\setminus\overline{\Omega}$ is Lipschitz and the Lipschitz character (number and Lipschitz constants of the local maps that cover the boundary) of $\Omdh$ is bounded independently of $h$.
  \end{enumerate}
\end{assumption}

The notation $a\eqsim b$ means that there exist constants $C,C'>0$ independent of $h$, but possibly dependent on $\kappa$ and $N$ such that $a\leq Cb$ and $b\leq C'a$, and similarly we write $a\lesssim b$ if and only if $a\leq Cb$.

\begin{remark}[Quasi-uniformity along cut-to-uncut paths]\label{rem:quasi.unif}
We remark that, by shape regularity and since $N$ is bounded independently of $h$, all mesh elements along a given cut-to-uncut path have uniformly comparable diameters: if $T_1,\ldots,T_{N_T}$ are as in \ref{item:cut-to-uncut}, then $h_{T_i}\eqsim h_{T_j}$ for all $i,j\in\{1,\ldots,N_T\}$.
\end{remark}
\section{Tangential and normal traces and parts of differential forms}\label{appendix:bdry_forms}

Restrictions as well as tangential and normal traces and parts of differential forms are essential concepts to define the stabilisation in \Cref{sec:ghost_penalty} below.
We consider here an element $T\in\mesh$ and one of its facets $F\subset\partial T$.

For $x\in F$, let $\Tan_x C$ be the tangent space of $C \in \{ F, T\}$ at $x$.
Let
\[
\iota:F\hookrightarrow T
\]
be the inclusion map of $F$ into $T$.
We note that $\Tan_xF$ is the vector hyperplane parallel to $F$ in $\RR^n$, and that $X\in \Tan_xF$ can be considered as a vector in $\RR^n$; for this reason, we identify $\iota_*X$ with $X$.

For $\omega\in \Lambda^kT$, by ``restriction to the facet $F$'' we simply mean the restriction of each component of $\omega$ to $F$, i.e., writing $\omega \coloneqq \sum_{1\leq j_1<j_2< \dots <j_k \leq n} \omega_{j_1 j_2 \dots j_k} \, dx^{j_1}\wedge dx^{j_2}\wedge \dots \wedge dx^{j_k}$, we define
\begin{equation}\label{eq:def.restriction}
    \omega|_F \coloneqq \sum_{1\leq j_1<j_2< \dots <j_k \leq n} (\omega_{j_1 j_2 \dots j_k}|_F) \, dx^{j_1}\wedge dx^{j_2}\wedge \dots \wedge dx^{j_k}.
\end{equation}
Notice that the result is a form which formally acts on $k$-tuples of vectors in $T$, although its components will only be evaluated on $F$. The reason is that we want $\omega|_F$ as a differential form to also be able to act on vectors normal to $F$.

The standard trace operator $\gamma \coloneqq \iota^* : \Lambda^kT \to \Lambda^kF$ on the facet is defined as the pullback of the inclusion $\iota$:
For all $\omega \in  \Lambda^k T$,%
\begin{equation}\label{eqs:trace}
  \forall x\in F,\qquad
  (\gamma \omega)_x(X_1,\dots ,X_k) \coloneqq \omega_x(X_1,\dots,X_k)
  \qquad \forall X_1,\dots,X_k\in \Tan_xF.
\end{equation}

Recall the Hodge star operator $\star:\Lambda^kT\to \Lambda^{n-k}T$.
We can also define a Hodge star $\star_F:\Lambda^kF\to \Lambda^{n-1-k}F$ on $F$ satisfying $\star_F\star_F \omega=(-1)^{k(n-1-k)}\omega$ for all $\omega\in \Lambda^kF$.

The normal trace operator is defined by ``sandwiching'' the standard trace between the Hodge star operators on $T$ and on $F$, see \cite[Eqs.~(45) and (46)]{bossavitDG}:
$\gamma_\normal \coloneqq (-1)^{n(k-1)} \star_F\gamma\star: \Lambda^kT\to \Lambda^{k-1}F$ is therefore such that, for all $\omega \in  \Lambda^k T$,
\begin{equation}\label{eqs:normaltrace}
  \forall x\in F,\qquad
  (\gamma_\normal \omega)_x(X_2,\dots ,X_k) \coloneqq
  \omega_x(\normal,X_2,\dots,X_k)
  \qquad \forall X_2,\dots,X_k\in \Tan_xF.
\end{equation}
The latter expression is, by definition, the contraction $\normal\lrcorner\, \omega$ of $ \omega$. 
The following identities hold, \cite[Eq.~(47)]{bossavitDG}:
\begin{align}\label{eqs:trace_relations}
  \star_F \gamma_\normal = \gamma \star,
  \qquad \gamma_\normal \star = (-1)^k\star_F\gamma.
\end{align}
From now on, we shall denote both $\star$ and $\star_F$ by $\star$, the context making it clear which one is meant.
The vector proxy interpretation of the trace and normal trace operators are summarised in Table~\ref{tab:trace.trace.n:proxy}.

\begin{table}[H]\centering
  \begin{tabular}{ccccc}
    \toprule
    Form degree $k$ & Proxy & Trace ($\gamma$) & Normal trace ($\gamma_\normal$) \\
    \midrule
    0 & $q$ & $q$ & 0
    \\
    1 & $\boldsymbol{v}$ & $\normal \times (\boldsymbol{v} \times \normal)$ & $\boldsymbol{v} \cdot \normal$
    \\
    2 & $\boldsymbol{w}$ & $\boldsymbol{w} \cdot \normal$ & $- \normal \times \boldsymbol{w}$
    \\
    3 & $r$ & $0$ & $r$ \\
    \bottomrule
  \end{tabular}
  \caption{Trace and normal trace for vector proxies in $\RR^3$.\label{tab:trace.trace.n:proxy}}
\end{table}

Sometimes, the standard trace is called the tangential trace, in contrast to the normal trace.
This nomenclature can be understood comparing the traces with the restriction for sufficiently smooth differential forms. 
Consider for the remainder of this section a smooth $k$-form $\omega\in\Lambda^kT$.
Using the pullback of the orthogonal projection $\pi_F: T\to F$, we can define the \emph{tangential part} of $\omega$ on $F$ by setting
\begin{equation}\label{eq:beta.parallel}
  \omega_\parallel \coloneqq (\pi_F^*\circ \gamma)\omega.
\end{equation}
At any $x \in F$, it is also defined as an alternating form in $\Alt^k \RR^n$.
Setting then
\begin{equation}\label{eq:beta.perp}
  \omega_\perp\coloneqq\omega|_F-\omega_\parallel,
\end{equation}
we have the following decomposition:
\begin{align}\label{eq:parallel_decomp}
  \omega|_F = \omega_\parallel + \omega_\perp.
\end{align}

We collect and prove some basic facts, see \cite[Eqs. (2.26)-(2.27)]{schwarz1995hodge}.
\begin{lemma}[Characterisation of the tangential and normal parts]\label{lem:normaltangential_components}
    Let $\omega\in \Lambda^kT$ be a smooth $k$-form. Then the following holds:
    \begin{align} \label{eq:beta.parallel=0<=>gamma.beta=0}
      \omega_\parallel = 0 &\iff \gamma\omega = 0,
      \\ \label{eq:beta.perp=0<=>gamma.n.beta=0}
     \omega_\perp = 0 &\iff \gamma_\normal\omega = 0.
    \end{align}
    Moreover,
    \begin{equation}\label{eq:star.traces}
      \text{%
        $(\star\omega)_\parallel = \star\omega_\perp \quad$ and
        $\quad (\star\omega)_\perp = \star\omega_\parallel$.
      }
    \end{equation}
\end{lemma}

\begin{proof}
  \underline{1) \emph{Proof of \eqref{eq:beta.parallel=0<=>gamma.beta=0}.}}
  We start by showing that $\omega_\parallel=0$ if and only if $\gamma \omega = 0$.
  It is clear by its definition~\eqref{eq:beta.parallel} that $\omega_\parallel=0$ if $\gamma\omega = 0$.
  To prove the converse, we just note that, if $(\omega_\parallel)_x=0$, then, for all $X_1,\dots,X_k\in \Tan_xF\subset \RR^n$, we have
  \[
  (\gamma\omega)_x(X_1,\dots,X_k)=(\gamma\omega)_{\pi_F(x)}(D\pi_F(X_1),\dots,D\pi_F(X_k))=
  (\omega_\parallel)_x(X_1,\dots,X_k)=0,
  \]
  the first equality following from the fact that $D\pi_F$ is the orthogonal projection on $\Tan_xF$.
  \smallskip\\
  \underline{2) \emph{Proof of \eqref{eq:star.traces}.}}
  Equalling the Hodge-star of the decomposition of $\omega$ and the decomposition of $\star\omega$ gives $\star\omega =\star\omega_\parallel+\star\omega_\perp = (\star\omega)_\parallel + (\star\omega)_\perp$. The relations \eqref{eq:star.traces} follow by applying vectors in $\RR^n$, respectively tangential and normal to $F$.
  \smallskip\\
  \underline{3) \emph{Proof of \eqref{eq:beta.perp=0<=>gamma.n.beta=0}.}}
  We simply write
  \[
    0 = \omega_\perp \iff 0 = \star\omega_\perp \overset{\eqref{eq:star.traces}}= (\star\omega)_\parallel \overset{\eqref{eq:beta.parallel=0<=>gamma.beta=0}}\iff 0 = \gamma(\star\omega) \overset{\eqref{eqs:trace_relations}}= \star \gamma_\normal \omega \iff 0=\gamma_\normal\omega,
  \]
  where the first and last equivalences are consequences of the bijectivity of the Hodge star.
\end{proof}

It can be helpful to look at the decomposition of a $k$-form $\omega\in\Lambda^kT$ in a coordinate system adapted to the facet $F$. Fix any $x\in F$.
Let $\{\normal,e_1,\dots,e_{n-1}\}$ be a basis of $\RR^n$ such that $\{e_1,\dots,e_{n-1}\}$ is a basis of $\Tan_xF$. Let the dual basis be $\{ dt,dx^1,\dots, dx^{n-1}\}$, where $dt$ is the dual to $\normal$.
Define the multi-index set
\[
\mcJ_\ell \coloneqq \left\{ (i_1,\dots, i_{\ell}) : 1\leq i_1 < i_2 < \dots < i_{\ell} \leq n-1 \right\}.
\]
Then, for any $x\in F$, we can write
\begin{align}\label{eq:normal_tangential_decomp_coord}
    (\omega)_x = (\omega_\perp)_x + (\omega_\parallel)_x = \sum_{J\in \mcJ_{k-1}} \omega_{t,J}(x) dt\wedge dx^J + \sum_{I\in \mcJ_{k}}\omega_I(x) dx^I,
\end{align}
where $\omega_{t,J}$ and $\omega_I$ are coefficients with respect to the chosen basis.
One can check (upon identifying $\gamma dx^{i}=dx^i$) that
\begin{equation}\label{eq:normal_trace_coord}
\text{%
  $(\gamma\omega)_x = \sum_{I\in \mcJ_{k}}\omega_I(x) dx^I$\quad  and  \quad
  $(\gamma_\normal\omega)_x = \sum_{J\in \mcJ_{k-1}} \omega_{t,J}(x) dx^J$.
  }
\end{equation}
These relations shed more light on the relations stated in \Cref{lem:normaltangential_components}. They are also useful for the next result.

\begin{lemma}\label{lem:normal_tangential_wedge}
    Let $\omega\in\Lambda^kT$ and $\eta\in\Lambda^{n-k}T$. The following relations hold at any $x\in F$:
    \[
      \omega_\parallel \wedge \eta_\parallel = 0,\qquad
      \omega_\perp \wedge \eta_\perp = 0.
    \]
\end{lemma}
\begin{proof}
    Writing the decomposition \eqref{eq:normal_tangential_decomp_coord} for both $\omega$ and $\eta$, it is clear that $\omega_\perp\wedge\eta_\perp$ vanishes due to $dt$ being present in both expansions.
    On the other hand, we have
    \begin{align*}
        \omega_\parallel \wedge \eta_\parallel &= \sum_{I\in\mcJ_k,\ L\in\mcJ_{n-k}} \omega_I\eta_L dx^I\wedge dx^L,
    \end{align*}
    which vanishes since one index at least is common to both $I$ and $L$, as can be seen by writing (identifying the indices in the families with the sets of these indices):
    \[
    |I\cap L| = |I|+|L|-|I\cup L| = k+(n-k)-|I\cup L| \geq n-(n-1)=1.\qedhere
        \]
\end{proof}


\section{Ghost penalty}\label{sec:ghost_penalty}

Using, in FEEC discretisations of weak formulations of PDEs,  the inner product $(\bullet,\bullet)_{\Omdh}$ on the active mesh $\Omdh$ in the discrete formulation is possible, but introduces a geometric error that is at least $O(h)$ due to the fact that the mesh is not fitted to the boundary $\partial\Omega$. Avoiding this geometric error requires to only rely on the inner product $(\bullet,\bullet)_\Omega$ in these formulations, which leads to ill-conditioned systems since DOFs attached to elements having a very small intersection with $\Omega$ have a very small impact on this inner product. The goal of the ghost penalty stabilisation term designed in this section is precisely to restore robustness while preserving optimal convergence rates.
The name ``ghost'' derives from the finite difference terminology of calling ``ghost cells'' the elements outside of the physical domain.
The main result of this section is \Cref{thm:equivalentnorms}, where we prove that the norm associated to the ghost term is equivalent to the $L^2(\Omdh)$-norm.

It should be mentioned that multiple equivalent approaches to stabilisation of CutFEMs have appeared over the years, see the review article \cite{burman2025cutfemreview}. 
To keep the presentation concise, we have picked the face based stabilisation of \cite[Remark 1]{Bu10Ghost}, which presents a particular challenge for the language of exterior calculus considering the unconventional traces involved. We shall see that normal and tangential parts of forms are precisely the right objects to consider inside ghost penalty face jumps.

\subsection{Jumps of forms}\label{sec:form_jumps}

Let $F=\overline T_1 \cap \overline T_2\in \mcFstab$ with normal $\normal$ pointing from $T_2$ into $T_1$.
Set, for $\omega\in \Vh{k}{r}$ and $i=1,2$,
\[
\omega_i \coloneqq  \omega|_{T_i},
\]
where $\omega|_{T_i}$ denotes the standard component-wise restriction of $\omega$ to $T_i$ on the canonical basis of $k$-alternating multi-linear maps, similar to \eqref{eq:def.restriction}. 
We define the jump over $F$ as
\begin{equation*} 
  \begin{aligned}
    [\cdot] : \Lambda^k(\overline T_1)\times \Lambda^k(\overline T_2) &\to C^\infty(F;\ \Alt^k\RR^n),\\
    \omega =(\omega_1,\omega_2)&\mapsto [\omega] \coloneqq  \omega_1|_F- \omega_2|_F.
  \end{aligned}
\end{equation*}

\begin{remark}[Traces and jumps]\label{rem:trace.jumps}
We can extend the normal and tangential traces \eqref{eqs:trace}--\eqref{eqs:normaltrace} fiberwise. Taking the tangential trace as an example, as $\gamma: C^\infty(F;\ \Alt^k\RR^n) \to \Lambda^k(F)$ by 
\[
(\gamma\omega)_x(X_1,\dots, X_k) = \omega_x(X_1,\dots,X_k)
\]
for $X_1,\dots,X_k\in T_xF$. 
Note also that this is not the identity, since $\Lambda^k(F) = C^{\infty}(F; \Alt^k F)$. As an example, take $F=\{x_1=0\}\subset \RR^2$ and $\omega = a\, dx^1 + b\, dx^2$. Then $\omega|_F = a|_F\, dx^1 + b|_F\, dx^2$ whereas $\gamma(\omega|_F) = b|_F\, dx^2.$
However, with this convention, we can see that $\gamma[\omega] = \gamma\omega_1 - \gamma\omega_2,$ since restriction is the identity on $F$. 
\end{remark}

From the linearity of the jump operator and the relation between the tangential part and the trace operator (cf. \Cref{lem:normaltangential_components}), it follows that only the normal part of the form (see \Cref{appendix:bdry_forms}) is involved in the jump.

\begin{prop}[Jump of discrete forms]\label{lem:discrete_jumps_characterization}
  For all $\omega\in \Vh{k}{r}$, it holds
  \[
    [\omega] = [\omega_{\perp}],
    \]
    where $\omega_\perp$ is the normal part of the form defined by~\eqref{eq:beta.perp}.
\end{prop}

\begin{proof}
  Let $F=\overline T_1\cap \overline T_2$.
  We have that $\omega_i \overset{\eqref{eq:parallel_decomp}}=\omega_{i\parallel}+\omega_{i\perp}$ on $F$ for $i=1,2$.
  The discrete $k$-forms are exactly defined by their degrees of freedom, so that $\gamma \omega$ is continuous across inter-element facets, i.e. $\gamma \omega_1-\gamma \omega_2=0$.
  Let $\{Y_1,Y_2,\dots, Y_k\}$ be an arbitrary set of $k$ vectors in $\RR^n$. By definition \eqref{eq:beta.parallel} of the tangential part of $\omega$, it holds
  \[
  (\omega_{i\parallel})_x(Y_1,\dots,Y_k)
  = (\gamma \omega_i)_x(D\pi_F Y_1,\dots,D\pi_F Y_k)
  \qquad \forall x\in F,
  \]
  where $\pi_F$ is the orthogonal projection onto $F$.
  By linearity of the jump operator, we have:
  \begin{equation}\label{eq:omega.par.zero}
  (\omega_{1\parallel}-\omega_{2\parallel})_x(Y_1,\dots,Y_k)
  = (\gamma \omega_1-\gamma \omega_2)_x(D\pi_F Y_1,\dots,D\pi_F Y_k)=0
  \qquad \forall x\in F.
  \end{equation}
  The discrete $k$-forms have degrees of freedom defined so that in particular $\gamma \omega$ is continuous across inter-element facets, i.e. $\gamma [\omega]=0$. 
  By \eqref{eq:beta.parallel=0<=>gamma.beta=0} and \eqref{eq:omega.par.zero}, we infer $[\omega]_\parallel = 0$, so that $[\omega] = [\omega]_\perp$. Since taking the normal component is linear and the same normal is used on both sides of $F$, we have $[\omega]_\perp = [\omega_\perp]$ and the statement follows.
\end{proof}

The result of Proposition~\ref{lem:discrete_jumps_characterization} is translated in terms of vector proxies in \Cref{tab:jump:proxy}.

\begin{table}\centering
  \begin{tabular}{cccc}
    \toprule
    $k$ & Proxy & Jump $[\omega_\perp]$ \\
    \midrule
    0 & $q$ & $0$
    \\
    1 & $\boldsymbol{v}$ & $[\boldsymbol{v} \cdot \normal]$
    \\
    2 & $\boldsymbol{w}$ & $[\normal \times (\boldsymbol{w} \times \normal)]$
    \\
    3 & $r$ & $[r]$ \\
    \bottomrule
  \end{tabular}
  \caption{Jumps $[\omega]=[\omega_\perp]$ for vector proxies in $\RR^3$.\label{tab:jump:proxy}}
\end{table}

\subsection{Ghost penalty stabilisation} 

For an integer $\ell\geq 0$, we define the $\ell$-th order directional derivative $\dirder : \Vh{k}{r}\to L^2 \Lambda^k \Omega_h$ along $\normal$ as
\begin{equation*} 
  \nabla^{(\ell)}_{\normal} \omega
  \coloneqq  \sum_{1\leq j_1<j_2< \dots <j_k \leq n} (\partial^{(\ell)}_{\normal} \omega_{j_1 j_2 \dots j_k}) \, dx^{j_1}\wedge dx^{j_2}\wedge \dots \wedge dx^{j_k},
\end{equation*}
where $\partial^{(\ell)}_{\normal} \omega_{j_1 j_2 \dots j_k}$ is the $\ell$-th order directional derivative of each component with respect to a chosen coordinate system.
The definition does not depend on the choice of coordinate system, and we also identify $\nabla_\normal^{0}\equiv \Id.$

For $\omega\in \Vh{k}{r}$, it is clear that $\nabla^{(\ell)}_\normal \omega$ is a polynomial form of degree $r-\ell$ and that $\dirder \omega=0$ for $\ell>r$.
Notice that the jump of the normal derivatives of a discrete $k$-form $\omega$ does not reduce, in general, to the jump of the normal part $[(\dirder \omega)_\perp]$, in contrast to \Cref{lem:discrete_jumps_characterization}.

Recall the definitions \eqref{eqs:trace} of trace and \eqref{eqs:normaltrace} of normal trace, as well as Remark \ref{rem:trace.jumps} on their applicability to jumps. Letting $h_F$ be the diameter of $F$, we can then define the (ghost penalty) stabilisation term $s : \Vh{k}{r} \times \Vh{k}{r} \to \RR$ as
\begin{equation}\label{eq:ghost_penalty}
  s(\omega,\zeta)
  \coloneqq  \sum_{F\in\mcFstab} \sum_{\ell=0}^r \mu h_F^{2\ell+1} \int_{F} \left( \gamma_\normal [\dirder \omega] \wedge \star \gamma_\normal [\dirder \zeta] + \gamma[\dirder \omega] \wedge \star\gamma [\dirder \zeta] \right),
\end{equation}
where $\mu>0$ is a user-defined penalty parameter.
The form $s$ is symmetric and bilinear.
Notice that $s$ can be extended to act on smooth enough $k$-forms (see \eqref{def:HmLamk} for the definition) and it holds, owing to the single-valuedness of traces,
\begin{equation} \label{eq:stab_vanishes_on_regular_forms}
  s(\zeta,\omega) = s(\omega,\zeta) = 0
  \qquad \forall (\zeta,\omega) \in H^{r+1}\Lambda^k \Omega_h \times \Vh{k}{r}.
\end{equation}
We define the ghost inner product and norm respectively as
\begin{equation}\label{eqs:def_stab_inner_product}
  \begin{alignedat}{2}
    (\omega,\zeta)_s &\coloneqq  (\omega,\zeta)_{\Omega} + s(\omega,\zeta) &\qquad& \forall (\omega,\zeta) \in \Vh{k}{r}\times \Vh{k}{r},
    \\
    \|\omega\|_s &\coloneqq  \sqrt{(\omega,\omega)_s} &\qquad& \forall \omega\in \Vh{k}{r}.
  \end{alignedat}
\end{equation}

\begin{remark}[Formula \eqref{eq:ghost_penalty} for vector proxies]\label{rem:FE.proxy.stab}
  Let $n=3$.
  The simplifications of the jumps for the case $\ell = 0$ given in \Cref{tab:jump:proxy} do not extend to the case $\ell \ge 1$.
  However, the directional derivative of a $k$-form is still a $k$-form, so, using the results of \Cref{tab:trace.trace.n:proxy}, we infer that the integrand in \eqref{eq:ghost_penalty} is in fact the inner product of the full jumps of the normal derivatives regardless of form order $k$, see \Cref{tab:normal_deriv:proxy}.

\begin{table}\centering
  \begin{tabular}{ccc}
    \toprule
    Form degree $k$ &
    \makecell{
      $s(\bullet,\bullet)$
    } \\
    \midrule
    0 or 3 & $\displaystyle\sum_{F\in\mcFstab}\sum_{\ell=0}^r \mu h_F^{2\ell+1} \int_F [\partial^{\ell}_\normal \bullet][\partial^{\ell}_\normal \bullet]$
    \\
    1 or 2 & $\displaystyle\sum_{F\in\mcFstab}\sum_{\ell=0}^r \mu h_F^{2\ell+1} \int_F [\partial^{\ell}_\normal \bullet ]\cdot [\partial^{\ell}_\normal \bullet ]$\\
    \bottomrule
  \end{tabular}
  \caption{Ghost penalty stabilisation \eqref{eq:ghost_penalty} for vector proxies in $\RR^3$ reduce to the standard ones, with vector valued jumps for $k=1,2$.\label{tab:normal_deriv:proxy}}
\end{table}

\end{remark}

The goal of the remainder of this subsection is to show the following equivalence result.

\begin{thm}[Uniform equivalence for $L^2$-norms]\label{thm:equivalentnorms}
  The norms $\|\bullet\|_{\Omdh}$ and $\|\bullet\|_s$ are equivalent on $\Vh{k}{r}$, uniformly in $h$.
\end{thm}

The following technical lemma is what informs the definition of the stabilisation term \eqref{eq:ghost_penalty}.

\begin{lemma}[Local control through ghost penalty]\label{lem:stabineq}
  Let $r \ge 1$ be an integer and $\omega \in \Vh{k}{r}$.
  Fix a boundary facet $F\in \mcFstab$ shared by the elements $T_1$ and $T_2$ of $\mcT_h$.
  Then, the following inequality holds:
  \begin{align}
    \|\omega\|^2_{T_1} &\lesssim \|\omega\|^2_{T_2} + \sum_{\ell=0}^r h_F^{2\ell+1} \int_{F}\left(\gamma_\normal [\dirder \omega] \wedge \star \gamma_\normal [\dirder \omega] + \gamma[\dirder \omega] \wedge \star\gamma [\dirder \omega] \right). \label{eq:fromT1toT2}
  \end{align}
\end{lemma}

\begin{proof}
  Consider the orthogonal projection onto the hyperplane containing $F$, which we also denote by $\pi_F$.
  First, we claim that we can assume to start with a simplex $T_1$ in which the vertex $\mathrm{s}_F$ opposite to $F$ is such that $\pi_F \mathrm{s}_F \in F$.
  Indeed, if $T_1$ is such that $\pi_F \mathrm{s}_F\not\in F$,  then consider the maximal simplex $T_1'\subseteq T_1$ with facet $F$ and opposite vertex $\mathrm{s}_F' \in T_1$ such that $\pi_F \mathrm{s}_F' \in F$.
  By a standard scaling argument as in \cite[Lemma 1.25]{di2020hho}, using norm equivalence for polynomial forms, we have that $\|\omega\|_{T_1}\lesssim \|\omega\|_{T_1'}$.

  Consider now a point $x\in T_1$, and its orthogonal projection $x_F$ on $F$.
  Let $t\in\RR$ be the distance from $x$ to $x_F$, such that $x=x_F+t\normal$. 
  For $i\in\{1,2\},$ since the restriction $w_i\coloneqq w|_{T_i}$ is a polynomial form, there is a canonical extension outside of $T_i$. 
  The Taylor expansion centered at $x_F$ of $\omega_1$, or of the extension of $\omega_2$ to $T_1$, is given by
  \begin{align*}
    (\omega_i)_x &= \sum_{\ell=0}^r \frac{(\dirder \omega_i)_{x_F}}{\ell!} t^\ell.
  \end{align*}
  Since each $\omega_i$ is a polynomial form of degree $r$, the Taylor expansion is exact.
  Taking the difference of the two expansions gives
  \[
  (\omega_1)_x
  = (\omega_2)_x + \sum_{\ell=0}^r \frac{(\dirder \omega_1)_{x_F}-(\dirder \omega_2)_{x_F}}{\ell!} t^\ell
  = (\omega_2)_x + \sum_{\ell=0}^r \frac{[\dirder \omega]_{x_F}}{\ell!} t^\ell.
  \]
  Taking the squared $L^2(T_1)$-norm and using triangle inequalities, we get
  \[
  \begin{aligned}
    \|\omega_1\|^2_{T_1} &\lesssim \left(\|\omega_2\|_{T_1} + \frac{|t|^{\ell}}{\ell!} \sum_{\ell=0}^r \|[\dirder \omega]_{x_F}\|_{T_1}\right)^2  \\
    &\lesssim \|\omega_2\|^2_{T_1} + h_F^{2\ell} \sum_{\ell=0}^r \|[\dirder \omega]_{x_F}\|^2_{T_1}  \\
    &\lesssim \|\omega_2\|^2_{T_2} + h_F^{2\ell} \sum_{\ell=0}^r \int_{T_1}[\dirder \omega]_{x_F}\wedge \star [\dirder \omega]_{x_F} ,
  \end{aligned}
  \]
  where, in the second inequality, we have used the Cauchy--Schwarz inequality on the norm together with $|t|\lesssim h_T\eqsim h_F$, and the conclusion follows from a polynomial transport argument (similar to the one used in the proof of \cite[Lemma 1.25]{di2020hho}) to replace $\|\omega_2\|^2_{T_1}$ with $\|\omega_2\|^2_{T_2}$, and the definition of the $L^2$-norm to expand the second term.

  Let $\xi$ be the $k$-form on $T_1$ such that, for all $x \in T_1$, $\xi_x\coloneq [\dirder \omega]_{x_F}$; as above, $x_F$ denotes the projection of $x$ onto $F$.
  From \eqref{eq:parallel_decomp} we can decompose $\xi$ into its normal and tangential parts inside $T_1$ as
  \[
    \xi = \xi_\perp + \xi_\parallel.
  \]
  By definition, $\xi_x=\xi_{x_F}$ for all $x\in T_1$.
  Define on $T_1$ the $(n-k)$-form $\eta\coloneq\star\xi$, for which we also have $\eta_x=\eta_{x_F}$ for all $x \in T_1$. 
  By Lemma \ref{lem:normaltangential_components}, we have $\star \xi_\parallel =\eta_\perp$ and $\star\xi_\perp=\eta_\parallel$ which, combined with \Cref{lem:normal_tangential_wedge}, gives 
  \[ \text{%
    $\xi_\perp \wedge \star \xi_\parallel = \xi_\perp \wedge \eta_\perp = 0$ and
    $\xi_\parallel \wedge \star \xi_\perp = \xi_\parallel \wedge \eta_\parallel = 0$.
  }
  \]
  Hence, at $x_F$ we can decompose the wedge product as
  \begin{equation*}
    \xi \wedge \eta = \xi_\perp \wedge \eta_\parallel +  \xi_\parallel\wedge\eta_\perp.
  \end{equation*}
  Decomposing $\xi$ using parallel and normal coordinates to $F$ as in \eqref{eq:normal_tangential_decomp_coord}, we write
  \[
    \xi_x = (\xi_\perp)_x + (\xi_\parallel)_x = \sum_{J\in\mcJ_{k-1}} \xi_{t,J}(x) dt\wedge dx^J + \sum_{I\in \mcJ_k}\xi_I(x) dx^I,
  \]
  and similarly for $\eta$.
  By Fubini's theorem, denoting by $\chi$ the characteristic function of $T_1$ and using a coordinate system $(t,y)$ with $t$ the coordinate along $\normal$ and $y$ the coordinate along the affine span of $F$, we get
  \begin{align*}
    \int_{T_1} \xi_\perp \wedge \eta_\parallel &= \int_{\RR^{n-1}}\int_{\RR} \chi(t,x_F) \sum_{J\in\mcJ_{k-1},I\in\mcJ_k}\xi_{t,J}(x_F)\eta_I(x_F) \ dt\wedge dx^J\wedge dx^I \\
    &= \int_{\RR^{n-1}} \sum_{J\in\mcJ_{k-1},I\in\mcJ_k} \xi_{t,J}(x_F)\eta_I(x_F) \left( \int_{\RR} \chi(t,x_F)  dt \right)  dx^J\wedge dx^I \\
    &= \int_F L_{x_F}\sum_{J\in\mcJ_{k-1},I\in\mcJ_k} \xi_{t,J}(x_F)\eta_I(x_F) dx^J\wedge dx^I \\
    &= \int_F L_{x_F} \sum_{J\in\mcJ_{k-1}} \xi_{t,J} dx^J \wedge \sum_{I\in\mcJ_k} \eta_I dx^I \\
    \overset{\eqref{eq:normal_trace_coord}}&= \int_F L_{x_F} \gamma_\normal \xi \wedge \gamma\eta,
  \end{align*}
  where, starting from the third equality, $L_{x_F}$ denotes the length of the segment at the vertical of $x_F$ in $T_1$.
  Performing similar computations for $\int_{T_1} \xi_\parallel\wedge\eta_\perp$ and accounting for the fact that $dx^I\wedge dt= (-1)^{k} dt\wedge dx^I$, we arrive at
  \begin{align*}
    \int_{T_1} \xi \wedge \eta &= \int_F L_{x_F}\ \gamma_\normal\xi \wedge \gamma\eta + (-1)^{k}\int_F L_{x_F}\ \gamma\xi \wedge \gamma_\normal\eta \\
    &= \int_F L_{x_F} \gamma_\normal [\dirder \omega] \wedge \gamma(\star[\dirder \omega]) + (-1)^{k}\int_F L_{x_F} \gamma[\dirder \omega] \wedge \gamma_\normal(\star[\dirder \omega]) \\
    \overset{\eqref{eqs:trace_relations}}&= \int_F L_{x_F} \gamma_\normal [\dirder \omega] \wedge \star \gamma_\normal [\dirder \omega] + \int_F L_{x_F} \gamma[\dirder \omega] \wedge \star\gamma [\dirder \omega] \\
    &\lesssim h_F \left(\int_F \gamma_\normal [\dirder \omega] \wedge \star \gamma_\normal [\dirder \omega] + \int_F \gamma[\dirder \omega] \wedge \star\gamma [\dirder \omega] \right),
  \end{align*}
  where we have expanded $\xi$ and $\eta$ according to the respective definitions in the second step while, in the conclusion, we have used the fact that $L_{x_F}\lesssim h_F$ and that the integrands are both positive multiples of the volume form on $F$.
\end{proof}

The next result shows that the ghost penalty term is enough to control the extension of $\omega$ to the active mesh. As soon as the local control of \Cref{lem:stabineq} is established, the global control follow standard arguments found in the known CutFEM literature, e.g. \cite[Lemma 6]{burman2012cutfemnitsche} or \cite[Lemma 2.20]{gurkan2019stabilized}. We provide a proof for completeness.

\begin{lemma}[Global control through ghost penalty]\label{corol:ghost_control}
  For all $\omega\in \Vh{k}{r}$, it holds
  \begin{equation*}
    \|\omega\|^2_{\Omdh} \lesssim
    \| \omega \|_s^2.
  \end{equation*}
\end{lemma}

\begin{proof}
  Let $\meshimm \coloneqq \mesh\setminus \meshcut$ be the set of immersed elements.
  We note that the norm on fully immersed elements is already controlled by the norm on $\Omega$:
  \begin{align*}
    \sum_{T\in\meshimm}\|\omega\|^2_{T} \leq \|\omega\|^2_\Omega.
  \end{align*}
  It remains to evaluate the norm on cut elements.
  \begin{figure}
    \centering
    \begin{tikzpicture}[scale=1.2, thick]
      \foreach \x in {0,1,2}
      \draw (\x,2) -- (\x,4);
      \foreach \y in {2,3,4}
      \draw (0,\y) -- (2,\y);
      \draw (0,2) -- (1,3);
      \draw (1,2) -- (2,3);
      \draw (0,4) -- (1,3);
      \draw (1,3) -- (2,4);

      \draw[red, thick]
      (-0.2,3.3) .. controls (0.8,3.0) and (1.7,3.4) .. (2.2,3.7);
      \node at (-0.4,3) {\dots};
      \node at (2.4,3) {\dots};
      \node at (1.0,1.5) {\vdots};
      \node at (1.3,3.75) {\(T_1\)};
      \node at (1.75,3.2) {\(T_2\)};
      \node at (1.3,2.65) {\(T_3\)};
    \end{tikzpicture}
    \caption{Both $T_1$ and $T_2$ are cut elements, but $T_2$ has a fully immersed neighbour $T_3$. The boundary $\partial\Omega$ is the red curve. Here $N=3$.}
    \label{fig:cut_elements}
  \end{figure}

  By Assumption \ref{assu:mesh}-(ii), starting from a cut element, the number of elements we must pass through to reach a fully immersed element is at most $N$, see \Cref{fig:cut_elements} for an illustration in the case $N=3$.
  Let $g : \meshcut \to \meshimm$ be the map which associates to each cut element $T \in \meshcut$ the first fully immersed element $g(T) \in \meshimm$ which is reached by passing through at most $N$ cut elements. For a given $T'\in\mesh$, any $T\in\mesh$ such that $T'=g(T)$ lies in a chain of elements of cardinality $\le N$; by Remark \ref{rem:quasi.unif}, all elements along this chain have diameters that are comparable to $h_{T'}$, and therefore lie in a ball of radius $\lesssim N h_{T'}\lesssim h_{T'}$. As this ball has measure $\eqsim h_{T'}^n$ and each $T$ in it has measure $\eqsim h_T^n\eqsim h_{T'}^n$, we infer that
  \begin{equation}\label{eq:estimate.gT}
  |g^{-1}(T')|\lesssim \frac{h_{T'}^n}{h_{T'}^n}= 1.
  \end{equation}

  By applying \Cref{lem:stabineq} at most $N$ times, we then have
  \begin{align*}
    \sum_{T\in\meshcut} \|\omega\|^2_{T} &\lesssim \sum_{T\in\meshcut} \|\omega\|^2_{g(T)} + \sum_{F\in\mcFstab} \sum_{\ell=0}^r h_F^{2\ell+1} \int_F\left(\gamma_\normal[ \dirder \omega]\wedge\star\gamma_\normal[ \dirder \omega] + \gamma[ \dirder \omega]\wedge\star\gamma[ \dirder \omega]\right) \\
    \overset{\eqref{eq:estimate.gT}}&\lesssim \sum_{T'\in\meshimm} \|\omega\|^2_{T'} + \sum_{F\in\mcFstab} \sum_{\ell=0}^r h_F^{2\ell+1} \int_F\left(\gamma_\normal[ \dirder \omega]\wedge\star\gamma_\normal[ \dirder \omega]+ \gamma[ \dirder \omega]\wedge\star \gamma[ \dirder \omega]\right).
  \end{align*}
  This completes the proof since $\sum_{T'\in\meshimm} \|\omega\|^2_{T'} \le \|\omega\|^2_{\Omega}$.
\end{proof}

\begin{proof}[Proof of \Cref{thm:equivalentnorms}]
  Choose an arbitrary $\omega \in \Vh{k}{r}$.
  In view of Lemma~\ref{corol:ghost_control}, we only have to show that $\| \omega \|_s \lesssim \|\omega\|_{\Omdh}$. That $\|\omega\|_\Omega \lesssim \|\omega\|_{\Omdh}$ is clear since $\Omega \subset \Omega_h$.
  Standard inverse and trace inequalities for polynomials grant
  \begin{equation*} 
  \|\dirder \omega\|^2_{F}
  \lesssim h_F^{-2\ell} \|\omega\|^2_{F}
  \lesssim h_F^{-(2\ell+1)} \|\omega\|^2_{T},
  \end{equation*}
  where $T$ is an element having $F$ as facet (note that $h_F \eqsim h_T$).
  Hence, using a triangle inequality followed by the above relation, the norm of the jumps can be estimated as follows:
  denoting by $T_1$ and $T_2$ the two elements on each side of $F$,
  \begin{align*}
    \|[\dirder \omega]\|_F^2 \lesssim \sum_{T\in \{T_1,T_2\}} h_F^{-(2\ell+1)} \|\omega\|^2_{T},
  \end{align*}
  from which $\| \omega \|_s \lesssim \|\omega\|_{\Omdh}$ follows.
\end{proof}

\begin{remark}[Macro stabilisation]\label{rem:macro_stab}
  Notice that the proof of \Cref{corol:ghost_control} does not require the stabilisation term to act on all facets in $\mcFstab$. It is enough to consider the facets that form a path from each cut element to a fully immersed element, as in the definition of the map $g$ in the proof of \Cref{corol:ghost_control}. Such a stabilisation term contains fewer summands, and leads naturally to the macro stabilisation procedure outlined in for example \cite[Section 6.1]{frachon2024stokes}.
\end{remark}

\subsection{Hodge decomposition}\label{sec:stab_mixed}

Recall the definition \eqref{eq:Hhk} of the space $\mathfrak{H}_h^k$ of discrete harmonic forms for FEEC, replacing $\Omega$ by $\Omdh$.
This space is constructed using orthogonality with respect to the standard inner product of $L^2 \Lambda^k \Omdh$.
Since we will use the ghost product $(\bullet,\bullet)_s$ (cf.~\eqref{eqs:def_stab_inner_product}) to define numerical schemes, we need to consistently modify the space of harmonic forms as follows:
\begin{align}
  \mathfrak{H}_s^k \coloneqq  \left\{\rho_h \in \Vh{k}{r} :\ d\rho_h=0, \ (\rho_h,d\tau_h)_s = 0 \qquad \forall \tau_h\in \Vh{k-1}{r+1}\right\}.
  \label{eq:ghost_harmonic_forms}
\end{align}
Since both $(\bullet,\bullet)_s$ and $(\bullet,\bullet)_{\Omdh}$ are inner products on $\Vh{k}{r}$, and both $\mathfrak{H}^k_s$ and $\mathfrak{H}_h^k$ are orthogonal to the same space, these spaces are isomorphic.
We also have
\begin{align*}
  \Ker d_h = d\Vh{k-1}{r+1} \operp_s \mathfrak{H}_s^k,
\end{align*}
where $\operp_s$ denotes the orthogonal complement with respect to the ghost product.

Correspondingly, we obtain the Hodge decomposition (compare with \eqref{eq:discrete_hodge_decomp}):
\begin{align}\label{eq:ghost_Hodge_decomp}
  \Vh{k}{r} = (\Ker d_h)^{\perp_s} \operp_s d\Vh{k-1}{r+1} \operp_s \mathfrak{H}_s^k.
\end{align}

\section{Further ingredients for error analysis}\label{sec:ingredients}
The solution to the continuous problem \eqref{eqs:mixed} is defined on the physical domain $\Omega$, while the unfitted discrete problem is defined on the active domain $\Omdh\supset\Omega$.
To compare discrete and continuous quantities one must extend the continuous solution to the active domain. 
The fully discrete error analysis introduced in \cite{di2018third}, also detailed in \cite[Appendix A]{di2020hho}, is a way to avoid error terms involving the extended continuous solution on the active domain. The idea is simply to split the error into an interpolation error and a discrete error, only write the latter in the equivalent $s$-norm, and then use the standard interpolation estimates to bound the former.

For $m\geq 0$ we define the spaces 
\begin{align}
  H^{m}\Lambda^k\Omega &\coloneqq \big\{ \alpha\in H\Lambda^k\Omega : \alpha_{j_1 \dots j_k} \in H^m(\Omega)\ \forall 1\le j_1 < \dots < j_k\le n \big\}, \label{def:HmLamk}\\
  H^{m,m}\Lambda^k\Omega &\coloneqq \big\{ \alpha\in H^m\Lambda^k\Omega : d\alpha \in H^m\Lambda^k\Omega \big\},
\end{align}
with associated graph norms $\|\cdot\|_{H^{m}\Omega}$ and $\|\cdot\|_{H^{m,m}\Omega}$.
For extension of the continuous solution we utilise the Sobolev--Stein extension operator $E: H^{m,m}\Lambda^k\Omega \to H^{m,m}\Lambda^k\Omdh$,
satisfying, for all $\alpha\in H^{m,m}\Lambda^k\Omega$, the following properties \cite[Eq.~(3.9),(3.16),Theorem~3.6]{hiptmair2012universalextension}:
\begin{subequations}\label{eqs:extension_op}
  \begin{align}
    E\alpha|_{\Omega} &= \alpha, \label{eq:extension_op_restr}\\
    \|E\alpha \|_{H^{m,m}\Omdh} &\lesssim \|\alpha\|_{H^{m,m}\Omega}.  \label{eq:extension_op_bound}
  \end{align}
\end{subequations}
    
\subsection{Interpolation estimates}
For the reader's convenience, we restate the interpolation estimates for the bounded cochain projection of FEEC, c.f. \cite[Sec. 7.5]{arnold_feec_2018}.
Let $|\bullet |_{H^{r}\Omdh}$ denote the $H^{r}$-seminorm on the active domain $\Omdh$. 

\begin{prop}(\cite[Theorem 7.6]{arnold_feec_2018})\label{prop:interpolation_estimates}
  Let $\alpha \in H^{r+1}\Lambda^k\Omdh$ for $r\geq 1$, and let $\Pih{k}{r}:H^{r+1}\Lambda^k\Omdh\to \Vh{k}{r}$ be the bounded cochain projection of \cite[Sec. 7.5]{arnold_feec_2018}. Then,
  \begin{equation}\label{eq:interpolation_estimates}
    \|\alpha - \Pih{k}{r}\alpha\|_{\Omdh} \lesssim h^{r+1} |\alpha|_{H^{r+1}\Omdh}.
  \end{equation}
\end{prop}
Since we have the inclusion $H^{m,m}\Lambda^k\Omdh\subset H^m\Lambda^k\Omdh$, these estimates hold also for forms in $H^{m,m}\Lambda^k\Omdh$.
Using the continuous extension \eqref{eqs:extension_op}, we can then interpolate a continuous form defined on $\Omega$ to a discrete form defined on the active domain $\Omdh$ (by first extending it, and by interpolating this extension), with norm bounded by the $H^{r+1}$-norm on $\Omega$.
\begin{cor}\label{cor:interpolation_estimates}
  Let $\alpha \in H^{r+1,r+1}\Lambda^k\Omega$ for $r\geq 1$. Then,
  $$
  \|E\alpha - \Pih{k}{r} E\alpha\|_{\Omdh} \lesssim h^{r+1}(\|\alpha\|_{H^{r+1}\Omega} + \|d\alpha\|_{H^{r+1}\Omega}).
  $$
\end{cor}
\begin{proof}
  We write
  \[
  \|E\alpha - \Pih{k}{r} E\alpha\|_{\Omdh}
  \overset{\eqref{eq:interpolation_estimates}}\lesssim
  h^{r+1} |E\alpha|_{H^{r+1}\Omdh}
  \overset{\eqref{eqs:extension_op}}\lesssim h^{r+1} \|\alpha\|_{H^{r+1,r+1}\Omega} \lesssim h^{r+1}(\|\alpha\|_{H^{r+1}\Omega} + \|d\alpha\|_{H^{r+1}\Omega}).\qedhere
  \]
\end{proof}

If using the trimmed spaces $\hatVh{k}{r}$, the same estimates hold with $\Pih{k}{r}$ replaced by the trimmed cochain projection of order $r$ acting on forms of the same $r+1$ regularity, and all with other appearances of $r$ replaced by $r-1$.

\subsection{Poincaré inequality on the active domain}\label{sec:poincare_on_comp_domain}
The Poincaré inequality has to be obtained on the active domain which is changing on every mesh size iteration. Thus for each $h$, the standard Poincaré constant $C_P=C_P(\Omdh)$ is different. A natural question is whether there is a Poincaré constant independent of $\Omdh$. An affirmative answer was given for any form order $k\ge 0$ in \cite{nilsson2026uniformly}, using a recursively defined cochain extension operator. 

\begin{thm}(\cite[Corollary 14]{nilsson2026uniformly})\label{lem:unif_discr_poincare}
  Assume that $\Omdh\setminus\overline{\Omega}$ is Lipschitz.
  There exists a constant $c_P>0$ independent of $h$ and $\Omdh$ such that, for any $q_h \in (\Ker d_h)^{\perp_s}\cap \Vh{k}{r}$,
  \begin{align*}
    \|q_h\|_s \leq c_P  \|d q_h\|_s.
  \end{align*}
\end{thm}

By \Cref{assu:mesh}-(iii), \Cref{lem:unif_discr_poincare} holds for the active domains $\Omdh$ of our unfitted meshes.

\section{CutFEEC scheme and convergence analysis}\label{sec:CutFEEC}
Consider once again the Hodge Laplace equation in mixed weak form \eqref{eqs:mixed}.
The unfitted method reads as follows:

Find $(\sigma_h,\eta_h,\lambda_h)\in \Vh{k-1}{r+1} \times \Vh{k}{r} \times \mathfrak{H}_s^k$ such that
\begin{equation}\label{eqs:unfitted_discrete_mixed}
  \begin{aligned}
    (\sigma_h, \tau_h)_s - (\eta_h, d\tau_h)_s &= 0  &&\qquad\forall \tau_h\in \Vh{k-1}{r+1},  \\
    (d \sigma_h, \zeta_h)_s + (d\eta_h,d\zeta_h)_s + (\lambda_h, \zeta_h)_s &= (f, \zeta_h)_\Omega  &&\qquad\forall \zeta_h\in \Vh{k}{r},  \\
    (\eta_h, \rho_h)_s &= 0  &&\qquad\forall \rho_h\in \mathfrak{H}_s^k.
  \end{aligned}
\end{equation}
We deliberately take the inner product involving $f$ over the physical domain $\Omega$ instead of the active domain $\Omdh$, to prevent the need for an unphysical extension of the data $f$.

\begin{remark}[On designing an unfitted scheme]
Our approach for the design of the discrete system can be summarised as follows. 
\begin{enumerate}
  \item Take, on the active mesh, a standard FEM discretisation with certain properties well suited to the continuous problem. 
  \item Add stabilisation to the $L^2$-inner products on the domain, which makes them equivalent to the $L^2$-norm on the active mesh (larger than the domain), while keeping intact any saddle point structure. 
  \item Inherit automatically the properties from the fitted setting to the unfitted setting.
\end{enumerate}
See the works \cite{frachon2024darcy}\cite{burman2024cutdivfree}\cite{lehrenfeld2025analysis}\cite{bouck2026tracefem}, listed chronologically according to the first public
version on arXiv, for examples of this approach.
\end{remark}

\begin{remark}[Divergence constraint for the Poisson equation in mixed form]\label{rem:poisson_div_constraint}
In the case $k = n$, the proposed method leads to a scheme which highlights some interesting features of the structure preserving viewpoint. It corresponds in this case to the Poisson equation in mixed form with homogeneous Dirichlet conditions:
Find $(\sigma_h,\eta_h)\in \Vh{n-1}{r+1} \times \Vh{n}{r}$ such that
\begin{alignat*}{2}
  (\sigma_h, \tau_h)_s - (\eta_h, \dive\tau_h)_s &= 0  &&\qquad\forall \tau_h\in \Vh{n-1}{r+1}, \\
  (\dive \sigma_h, \zeta_h)_s &= (f, \zeta_h)_\Omega  &&\qquad\forall \zeta_h\in \Vh{n}{r}.
\end{alignat*}
Up to interface conditions and a source term, this is the formulation given in \cite[(3.10)-(3.11)]{frachon2024darcy}. 
If $f\in L^2\Lambda^n\Omega$ has $r+1$ weak normal derivatives in $L^2$, so that $f$ admits a sufficiently regular extension $F$ for which the ghost stabilisation vanishes
against discrete test functions, then $\dive\sigma_h$ will be equal, to machine precision, to the $L^2$-projection of $F$ with respect to the $(\bullet,\bullet)_s$-inner product.
Applying the same principle for Stokes flow ($f=0$), the corresponding scheme satisfies the incompressibility condition exactly since $\dive \sigma_h=0$. 
This observation, which is the content of \cite[Theorem~1]{frachon2024stokes}, can also be seen as a consequence of \Cref{thm:equivalentnorms}.
\end{remark}

We define:
\[
X_h \coloneqq  \Vh{k-1}{r+1}\times \Vh{k}{r}\times \mathfrak{H}_s^k,\qquad
X \coloneqq  H\Lambda^{k-1}\Omdh\times H\Lambda^{k}\Omdh\times \mathfrak{H}^{k}.
\]
Let $U_h \coloneqq  (\sigma_h,\eta_h,\lambda_h)\in X_h$ and $W_h \coloneqq  (\tau_h,\zeta_h,\rho_h)\in X_h$.
Define the bilinear form $\mcA_s: X_h \times X_h \to \RR$ and the linear form $\mcL_h(W_h) :X_h\to\RR$ such that, for all $(U_h,W_h) \in X_h \times X_h$,
\begin{alignat*}{2}
  \mcA_s(U_h,W_h)
  &\coloneqq  (\sigma_h, \tau_h)_s - (\eta_h, d\tau_h)_s + (d \sigma_h, \zeta_h)_s
  + (d\eta_h,d\zeta_h)_s + (\lambda_h, \zeta_h)_s- (\eta_h, \rho_h)_s,\\
  \mcL_h(W_h) &\coloneqq  (f, \zeta_h)_\Omega.
\end{alignat*}
The problem \eqref{eqs:unfitted_discrete_mixed} can be recast as: find $U_h\in X_h$ such that
\[
\mathcal{A}_s(U_h, W_h) = \mcL_h(W_h)\quad \forall W_h\in X_h.
\]
We also define the combined norms
\begin{align*}
  \|U_h\|^2_{\Omdh} &\coloneqq  \|\sigma_h\|^2_{\Omdh} + \|\eta_h\|^2_{\Omdh} + \|\lambda_h\|^2_{\Omdh}, \\
  \|U_h\|^2_{s} &\coloneqq  \|\sigma_h\|^2_{s} + \|\eta_h\|^2_{s} + \|\lambda_h\|^2_{s}, \\
  \|U_h\|^2_{d,s} &\coloneq \|U_h\|^2_s + \|dU_h\|^2_s.
\end{align*}


\subsection{Stability}\label{sec:stability} 

The stability of the scheme classically follows from an inf-sup condition for the bilinear form $\mcA_s$. To prove it, we follow the approach in \cite[Theorem 4.9]{arnold_feec_2018}, the main difference being the use of the ghost product.

\begin{lemma}[Inf-sup condition]\label{lem:infsup}
  For any $U_h\in X_h$ there exists $\check W_h\in X_h$ such that,
  \begin{align*}
    \frac{\mcA_s(U_h,\check W_h)}{\|\check W_h\|_{d,s}} \gtrsim \|U_h\|_{d,s}
  \end{align*}
  where the hidden constant contains the maximum of the Poincaré constants for discrete $k$- and $(k-1)$-forms.
\end{lemma}

\begin{proof}
 Using \eqref{eq:ghost_Hodge_decomp}, we write the $s$-orthogonal Hodge decomposition $\eta_h = \tilde \eta_h + d\omega_h + \xi_h$ with $\tilde \eta_h\in(\Ker d^k_h)^{\perp_s}$, $\xi_h\in\mathfrak{H}_s^k$ and $\omega_h\in (\Ker d_h^{k-1})^{\perp_s}$.
  Choosing $W_h = c_P^2U_h$ gives by the Poincaré inequality (\Cref{lem:unif_discr_poincare}) on $\tilde\eta_h$ and noticing that $d\eta_h=d\tilde\eta_h$:
  \begin{align*}
    \mcA_s(U_h,c_P^2U_h) &= c_P^2\|\sigma_h\|^2_{s} + c_P^2\|d\eta_h\|^2_{s} \\
    &= c_P^2\|\sigma_h\|^2_{s} + \frac{c_P^2}{2}\|d\eta_h\|^2_{s} + \frac{c_P^2}{2}\|d\eta_h\|^2_{s} \\
    &\ge c_P^2\|\sigma_h\|^2_{s} + \frac{c_P^2}{2}\|d\eta_h\|^2_{s} + \frac{1}{2}\|\tilde \eta_h\|^2_{s}.
  \end{align*}
  Choosing $W_h = (0,d\sigma_h+\lambda_h,0)$ gives, by $dd\sigma=0$ and $s$-harmonicity of $\lambda_h$ (which entails $(d\sigma_h,\lambda_h)_s=0$),
  \begin{align*}
    \mcA_s(U_h,W_h) &= \|d\sigma_h\|^2_s + \|\lambda_h\|^2_s.
  \end{align*}
  Lastly, choosing $W_h = (-\omega_h,0,-\xi_h)$ gives by $s$-orthogonality of the Hodge decomposition:
  \[
  \mcA_s(U_h,W_h) = -(\sigma_h, \omega_h)_s + (\eta_h, d\omega_h)_s -(\eta_h,-\xi_h)_s = -(\sigma_h, \omega_h)_s + \| d\omega_h\|^2_s + \|\xi_h\|^2_s.
  \]
  By Young's inequality and then the Poincaré inequality (\Cref{lem:unif_discr_poincare}) on $\omega_h\in (\Ker d^{k-1}_h)^{\perp_s}$,
  \[
  -(\sigma_h, \omega_h)_s
  \ge -\frac{c^2_P}{2}\|\sigma_h\|^2_s - \frac{1}{2c^2_P} \|\omega_h\|^2_s
  \ge -\frac{c^2_P}{2}\|\sigma_h\|^2_s - \frac{1}{2}\|d\omega_h\|^2_s
  \]
  so that
  \[
   \mcA_s(U_h,W_h) \ge -\frac{c^2_P}{2}\|\sigma_h\|^2_s  + \frac12\| d\omega_h\|^2_s + \|\xi_h\|^2_s.
  \]
  Hence, choosing $\check W_h = (c^2_P\sigma_h-\omega_h,c^2_P \eta_h+d\sigma_h+\lambda_h,c^2_P\lambda_h-\xi_h)$ as the sum of all previous choices gives a lower bound corresponding to the sum of all previous bounds:
  \begin{align*}
    \mcA_s(U_h,\check W_h)
    &\ge \frac{c_P^2}{2}\|\sigma_h\|^2_{s}+\|d\sigma_h\|^2_{s} + \frac12\|\tilde \eta_h\|^2_{s} + \|d\omega_h\|^2_s + \|\xi_h\|^2_s + \frac{c_P^2}{2}\|d\eta_h\|^2_s + \|\lambda_h\|^2_{s} \\
    &\gtrsim \|\sigma_h\|^2_{s} + \|d\sigma_h\|^2_{s} + \|\eta_h\|^2_{s} + \|d\eta_h\|^2_{s} + \|\lambda_h\|^2_{s} = \|U_h\|^2_s + \|dU_h\|^2_s.
  \end{align*}
  The proof is complete by noticing that $\|\check W_h \|_s + \|d\check W_h\|_s \lesssim (1+c_P^2)(\|U_h\|_s+\|dU_h\|_s)$.
\end{proof}


\subsection{Consistency error}\label{sec:consistency}
Let $U\in X = H\Lambda^{k-1}\Omega \times H\Lambda^k\Omega \times \mathfrak{H}^k$ be the solution to the Hodge Laplace equation \eqref{eqs:mixed}.
Let $\hat U \coloneq EU$ be the extension of $U$ to the active domain, where $E$ is the Sobolev--Stein extension operator \eqref{eqs:extension_op} applied to each component of $U$.
Define also the projection operator 
$$
  \Pi_h\coloneqq (\Pih{k-1}{r+1},\Pih{k}{r},\Pih{k}{r}): X\to X_h.
$$
For the sake of notational simplicity, we will in the following not specify the order $r$ of the cochain projections, and just write $\Pi_h^k$ instead of $\Pih{k}{r}$, and similarly for $\Pi_h^{k-1}$. 
The consistency error in the fully discrete framework is defined by the linear form $\mcE_h(U,\bullet): X_h \to \RR$ defined by
\begin{align*}
  \mcE_h(U,W_h) &\coloneqq  \mcL_h(W_h) - \mcA_s(\Pi_h \hat U,W_h).
\end{align*}


By the third equation of \eqref{eqs:mixed}, $\eta$ is $L^2$-orthogonal to harmonic forms on $\Omega$, so that by the continuous Hodge decomposition, we have $\eta = \tilde \eta + d\omega$, where $\tilde \eta\in (\Ker d|_{\Omega})^{\perp_\Omega}$ and $\omega\in H\Lambda^{k-1}\Omega$.

\begin{lemma}\label{lem:orthogonality}
  Assume the solution $\tilde \eta+d\omega\in H\Lambda^k\Omega$ to the Hodge Laplace equation \eqref{eqs:HodgeLaplace}
  is such that $\tilde \eta$ and $\omega$ are both of regularity $H^{r+1,r+1}(\Omega)$. Then
  \begin{align}
  \label{eq:Eeta.perp}
  E\tilde \eta &\perp_s \mathfrak{H}^k_s,\\
  \label{eq:dEeta.perp}
  \Pi_h^k dE\omega &\perp_s \mathfrak{H}^k_s.
  \end{align}
\end{lemma}
\begin{proof}
  Thanks to the regularity of $\eta$, the extension $E\eta$ belongs to $H^{r+1,r+1}\Lambda^k\Omdh$ so that $s(E\tilde \eta,\bullet)_s=0$ by \eqref{eq:stab_vanishes_on_regular_forms}.
  Let $\rho_h\in \mathfrak{H}^k_s$ and notice that $d\rho_h = 0$ by definition \eqref{eq:ghost_harmonic_forms}, implying $\rho_h\in\Ker d|_\Omega$. 
  Then, since $\tilde \eta\in (\Ker d|_{\Omega})^{\perp_\Omega}$, we obtain
  $(E\tilde \eta, \rho_h)_s = (\tilde \eta, \rho_h)_\Omega = 0$,
  which shows that $E\tilde \eta \perp_s \mathfrak{H}^k_s$.

  By the commuting property of the cochain projection, we have $\Pi_h^k dE\omega = d\Pi_h^{k-1} E\omega \in d\Vh{k-1}{r+1}$. Therefore $\Pi_h^k dE\omega \perp_s \mathfrak{H}^k_s$ by the ghostly Hodge decomposition \eqref{eq:ghost_Hodge_decomp}.
\end{proof}

\Cref{lem:orthogonality} gives us the extension we must consider for the solution $\eta$ of \eqref{eqs:mixed}:
\begin{align*}
  \hat \eta &\coloneqq E\tilde \eta + dE\omega.
\end{align*}
Note that $\hat\eta$ is not equal to $E\eta$ in general, since $E$ does not commute with $d$. Moreover, the following equation is satisfied for any $\rho_h\in \mathfrak{H}^k_s$ thanks to \eqref{eq:Eeta.perp} and \eqref{eq:dEeta.perp}:
\begin{align}
  (\Pi^k_h\hat \eta,\rho_h)_s &= (\Pi^k_h E\tilde \eta - E\tilde \eta,\rho_h)_s.\label{eq:nec_consistency_diff}
\end{align}
Let also $\widehat{\sigma} \coloneqq E\sigma$, $\widehat{\lambda} \coloneqq E\lambda,$ and set 
$
  \widehat{U}
  \coloneqq
  \left(
      \widehat{\sigma},
      \widehat{\eta},
      \widehat{\lambda}
  \right).
$
\begin{lemma}[Consistency estimate]\label{lem:consistency}
  Assume the solution $U\in X$ to the Hodge Laplace equation \eqref{eqs:mixed} is of regularity $H^{r+1,r+1}(\Omega)$, and assume further in the Hodge decomposition of $\eta$ that $\tilde \eta$ and $\omega$ are of regularity $H^{r+1,r+1}(\Omega)$.
  Then for all $W_h\in X_h$,
  $$ \mcE_h(U,W_h) \lesssim h^{r+1} (\|\sigma\|_{H^{r+1}\Omega}+\|d\sigma\|_{H^{r+1}\Omega}+\|\eta\|_{H^{r+1}\Omega}+\|d\eta\|_{H^{r+1}\Omega}+\|\lambda\|_{H^{r+1}\Omega})(\|W_h\|_s + \|dW_h\|_s) .$$
\end{lemma}

\begin{proof}
  Since $\Pi_h^k$ is a cochain map we have,
  \begin{align}
    \mcA_s(\Pi_h \hat U,W_h)
    &= (\Pi_h^{k-1} \hat\sigma, \tau_h)_s - (\Pi_h^k \hat \eta, d\tau_h)_s + (\Pi_h^k d \hat \sigma, \zeta_h)_s
    + (\Pi_h^{k+1} d \hat \eta,d\zeta_h)_s + (\Pi_h^k \hat \lambda, \zeta_h)_s- (\Pi_h^k \hat \eta, \rho_h)_s. \label{eq:mcAs_expression}
  \end{align}
  The regularity of $U$ implies by \eqref{eq:stab_vanishes_on_regular_forms}, that the stabilisation terms vanish:
  \begin{align}\label{eq:stabilisation_vanishes}
    s(\hat \sigma,\bullet)_s = s(\hat \eta,\bullet)_s=s(d\hat \sigma,\bullet)_s=s(d\hat \eta,\bullet)_s =s(\hat\lambda,\bullet)=0.
  \end{align}
  By \eqref{eq:nec_consistency_diff}, which holds thanks to the regularity assumptions, we have that
  \begin{align*}
    -(\Pi_h^k \hat \eta, \rho_h)_s &= (E\tilde \eta - \Pi_h^k E\tilde \eta, \rho_h)_s.
  \end{align*}
  We now use that the continuous solution satisfies the equation $f=\pi f + d\delta \eta + \delta d\eta = \lambda + d\sigma + \delta d\eta$. Integrating by parts over the physical domain $\Omega$ and using the second boundary condition $\gamma \star d\eta = 0$, we have
  \begin{align*}
    \mcL_h(W_h) &= (f,\zeta_h)_\Omega = (d\sigma,\zeta_h)_\Omega + (d\eta,d\zeta_h)_\Omega + (\lambda, \zeta_h)_\Omega.
  \end{align*}
  We can also integrate $\sigma- \delta \eta=0$ by parts along with the first boundary condition $\gamma \star \eta = 0$ to get,
  \begin{align*}
    \mcL_h(W_h) &= (\sigma,\tau_h)_\Omega - (\eta,d\tau_h)_\Omega + (d\sigma,\zeta_h)_\Omega + (d\eta,d\zeta_h)_\Omega + (\lambda,\zeta_h)_\Omega \\
    &= ( \hat \sigma,\tau_h)_s - ( \hat \eta,d\tau_h)_s + (d \hat \sigma,\zeta_h)_s + (d \hat \eta,d\zeta_h)_s + ( \hat \lambda,\zeta_h)_s, 
  \end{align*}
  the last equality being due \eqref{eq:stabilisation_vanishes}.
  Thus, combining with \eqref{eq:mcAs_expression}, we get 
  \begin{align*}
    \mcE_h(U,W_h) &= (\hat \sigma - \Pi_h^{k-1} \hat\sigma, \tau_h)_s - (\hat \eta - \Pi_h^k \hat \eta, d\tau_h)_s + (d\hat \sigma - \Pi_h^k d\hat \sigma,\zeta_h)_s\\
    &\qquad + (d\hat \eta - \Pi_h^{k+1} d\hat \eta,d\zeta_h)_s + (\hat \lambda - \Pi_h^k \hat \lambda, \zeta_h)_s 
    + (E\tilde \eta - \Pi_h^k E\tilde \eta, \rho_h)_s \\
    &\lesssim (\|\hat \sigma - \Pi_h^{k-1} \hat\sigma\|_s + \|\hat \eta - \Pi_h^k \hat \eta\|_s + \|d\hat \sigma - \Pi_h^k d\hat \sigma\|_s + \|d\hat \eta - \Pi_h^{k+1} d\hat \eta\|_s + \|\hat \lambda - \Pi_h^k \hat \lambda\|_s
    )\\
    &\qquad\times
    (\|\tau_h\|_s+\|d\tau_h\|_s+\|\zeta_h\|_s+\|d\zeta_h\|_s+\|\rho_h\|_s).
  \end{align*}
  The interpolation estimates \Cref{cor:interpolation_estimates} together with the norms equivalence of \Cref{thm:equivalentnorms} concludes the proof.
\end{proof}

\subsection{Convergence by Third Strang Lemma}\label{sec:convergence} 

Using the triangle inequality, split the error into an interpolation error and a discrete error and use Lemma \ref{corol:ghost_control} to bound the $L^2$-norms of discrete forms by the stabilised norm:
\begin{align*}
  \|\hat U - U_h\|_{H\Omdh} \lesssim  \|\hat U-\Pi_h \hat U\|_{H\Omdh}+\|U_h-\Pi_h \hat U\|_s+\|dU_h-\Pi_hd \hat U\|_s.
\end{align*}
The convergence of the first term follows from standard interpolation estimates. The last two terms together form a discrete $H\Lambda^k$-norm difference on the active mesh, which can be bounded using \Cref{lem:infsup} and \Cref{lem:consistency} and invoking the Third Strang lemma \cite[Theorem 10]{di2018third}. To summarise, we obtain the following theorem.

\begin{thm}
  Under the assumptions of Lemma \ref{lem:consistency}, the solution $U_h\in X_h$ to the mixed method \eqref{eqs:unfitted_discrete_mixed} converges in $H\Omega$-norm to $U$ with order $O(h^{r+1})$.
\end{thm}

\begin{remark}[Required properties of stabilisation]
We highlight three main properties inherent to the stabilisation operator \eqref{eqs:def_stab_inner_product} that were sufficient to perform the error analysis in this section. The idea is that any other stabilisation operator satisfying the following three properties would be equally well suited to discretise the equations.
\begin{enumerate}
  \item \Cref{thm:equivalentnorms}. It should when added to the standard physical domain $L^2$-norm constitute an equivalent norm to the $L^2$-norm over the active mesh. This allows also to obtain the Hodge decomposition \eqref{eq:ghost_Hodge_decomp} with respect to the stabilised inner product.
  \item \Cref{lem:unif_discr_poincare}. Necessary for stability is the uniform Poincaré inequality.
  \item \Cref{eq:stab_vanishes_on_regular_forms}. To perform the a priori estimate, the continuous forms need to lie in the kernel of the stabilisation operator.
\end{enumerate}
\end{remark}

\section{Numerical validation}\label{sec:numerics}
In this section we validate our theory with two examples aimed to illustrate the utility of the structure preserving approach. The first example is a convergence study for the Hodge Laplace equation on a filled torus, which has a nontrivial harmonic form. We also show that the condition number of the system matrix is independent of the position of the domain with respect to the background mesh. 
The second example considers a simple $H^{\crl}$-conforming method for Stokes which is known to be pressure robust in the fitted setting \cite{da2022presrobstokes}, and we show that the same property holds in the unfitted setting.

The C++ code used to compute the example can be found in the fork \cite{tensorFa11:online} of the CutFEM library \cite{CutFEMCu93:online}. Linear systems were solved using the sparse direct solver MUMPS \cite{amestoy2001mumps}.

\subsection{Implementation and results}

We validate in this section the scheme \eqref{eqs:unfitted_discrete_mixed} for $k=1$.
We respectively choose the trimmed spaces $P^-_{1}\Lambda^{0}\Omdh$ (Lagrange) and $P^-_1\Lambda^1\Omdh$ (Nédélec of the first kind) for $\sigma_h$ and $\eta_h$.
This choice means that the exterior derivative $d$ does not reduce the indicated polynomial order, and that the error in the natural stability norm can be expected to converge with order $O(h)$. 
The scheme reads as follows.
Find $(\sigma_h,\eta_h,\lambda_h)\in P^-_{1}\Lambda^{0}\Omdh \times P^-_1\Lambda^1\Omdh \times \mathfrak{H}_s^1$ such that
\begin{equation}\label{eq:discrete:k=1}
  \begin{alignedat}{2}
    (\sigma_h, \tau_h)_s - (\eta_h, \nabla\tau_h)_s &= 0  &&\qquad\forall \tau_h\in P^-_1\Lambda^0\Omdh,  \\
    (\nabla \sigma_h, \zeta_h)_s + (\crl \eta_h,\crl \zeta_h)_s + (\lambda_h, \zeta_h)_s &= (f, \zeta_h)_\Omega  &&\qquad\forall \zeta_h\in P^-_1\Lambda^1\Omdh,  \\
    (\eta_h, \rho_h)_s &= 0  &&\qquad\forall \rho_h\in \mathfrak{H}_s^1,
  \end{alignedat}
\end{equation}
where we use stabilisation parameter $\mu=1$ in \eqref{eq:ghost_penalty}.
Since the trial and test forms are piecewise linear, their second derivatives vanish identically, and since $\sigma_h,\tau_h$ are continuous, we get from Tables~\ref{tab:jump:proxy} and \ref{tab:normal_deriv:proxy} that 
\begin{align*}
  (\sigma_h, \tau_h)_s &= (\sigma_h, \tau_h)_\Omega + \sum_{F\in \mcFstab}h_F^3([\nabla\sigma_h\cdot\normal],[\nabla\tau_h\cdot\normal])_{F}, \\
  (\eta_h,\nabla\tau_h)_s &= (\eta_h,\nabla\tau_h)_\Omega + \sum_{F\in \mcFstab}h_F([\eta_h],[\nabla\tau_h])_{F},\\
  (\crl \eta_h,\crl \zeta_h)_s &= (\crl \eta_h,\crl \zeta_h)_\Omega + \sum_{F\in \mcFstab}h_F([\crl \eta_h],[\crl \zeta_h])_{F}, \\
  (\lambda_h, \zeta_h)_s &= (\lambda_h, \zeta_h)_\Omega + \sum_{F\in \mcFstab}\left( h_F([\lambda_h], [\zeta_h])_{F} + h_F^3([\nabla\lambda_h\cdot\normal], [\nabla \zeta_h\cdot\normal])_{F} \right).
\end{align*}
We also use the macro parameter $\delta=0.25$ in the macro stabilisation procedure discussed in \Cref{rem:macro_stab}.

\begin{figure}
  \centering
  \includegraphics[width=0.95\textwidth]{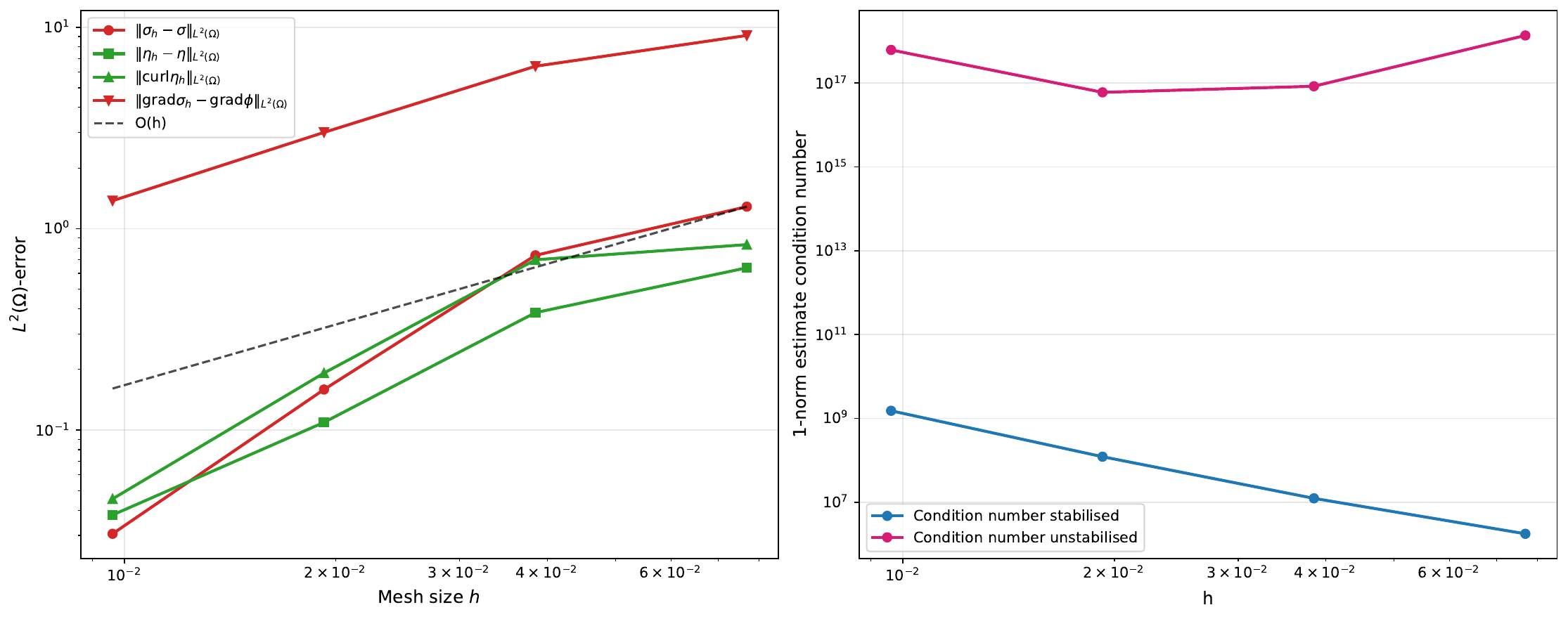}
  \caption{
    Hodge Laplace example: Plot of the convergence of the mixed method \eqref{eqs:unfitted_discrete_mixed} as a function of the mesh size $h$. The convergence is shown for polynomial order $r=1$, and the $1$-norm estimate of the condition number $\kappa_s \coloneq \|\mcA_s\|_{\op} \|\mcA_s^{-1}\|_{\op}$ of the system matrix $\mcA_s$ is also plotted on the right figure.
  }
  \label{fig:error_and_condition_plots}
\end{figure}

\begin{figure}
  \centering
  \includegraphics[width=0.75\textwidth]{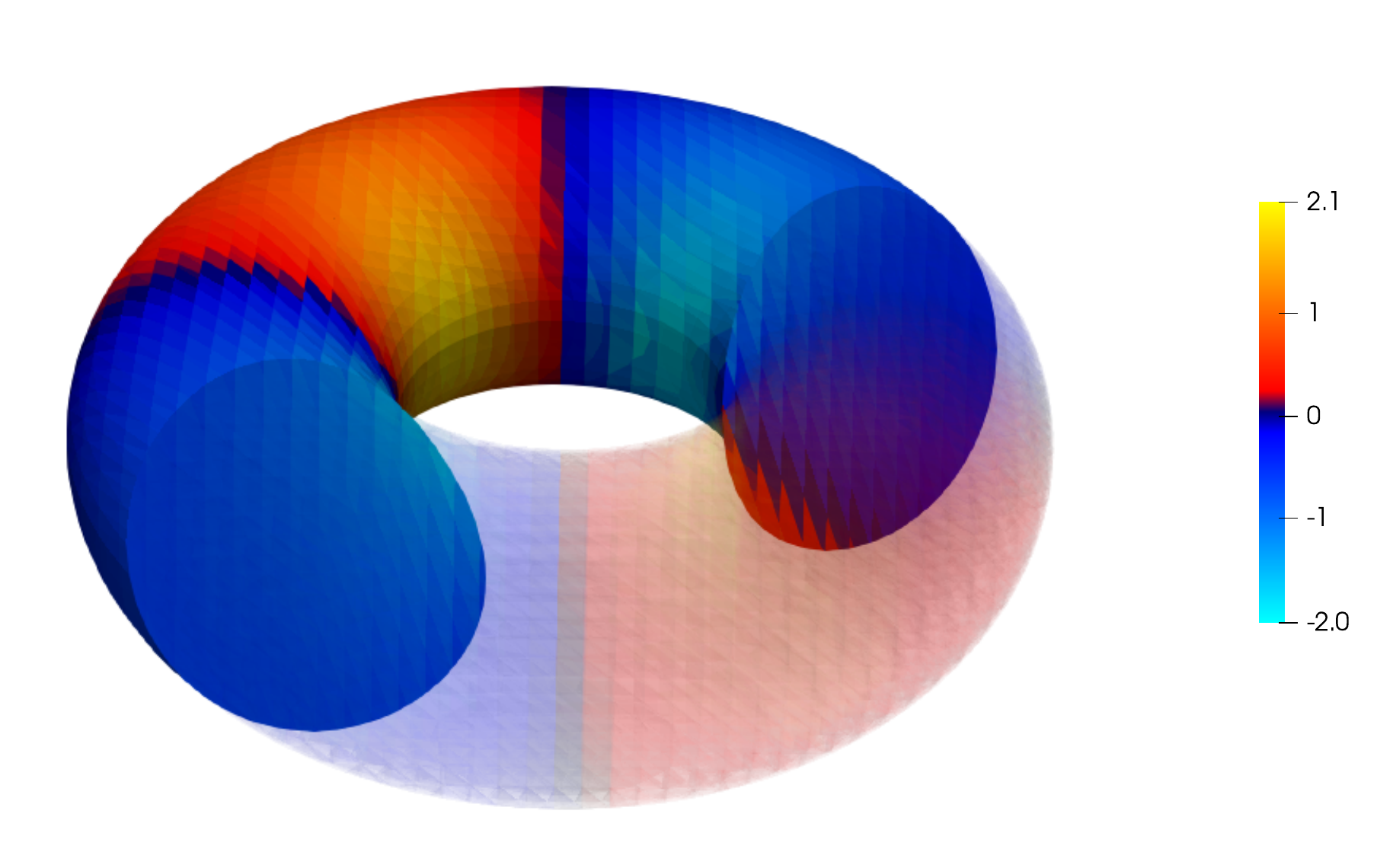}
  \caption{
    Hodge Laplace example: Heat map of computed $x$-component of the solution $\eta_h$ to the discrete Hodge Laplace equation \eqref{eqs:unfitted_discrete_mixed} on the filled in torus. The mesh size is $h=0.019231$ and the polynomial order is $r=1$. 
  }
  \label{fig:HL_ux}
\end{figure}

We solve problem~\eqref{eq:discrete:k=1} on the following filled torus embedded in $\RR^3$: 
$\Omega = \{(x,y,z)\in\RR^3: [(\sqrt{x^2+y^2}-0.5)^2 + z^2]^{1/2}\leq 0.25\}$, with major radius $0.5$ and minor radius $0.25$.
The mesh is cut via the level set function
\begin{align*}
  \phi = \left[ \left(\sqrt{x^2+y^2}-0.5\right)^2 + z^2\right]^{1/2} - 0.25,
\end{align*}
from a uniform background mesh with mesh size $h$.
The toroidal harmonic forms are scalar multiples of
\[
  d\theta = \frac{xdy-ydx}{x^2+y^2} \longleftrightarrow \frac{1}{x^2+y^2}(-y,x,0).
\]
We pick the right-hand side and the exact solution as
\[
f= \left(-\frac{3xy}{(x^2+y^2)^{\frac52}}, \frac{x^2-2y^2}{(x^2+y^2)^{\frac52}}, 0\right), \quad\quad \eta = \left(-\frac{-xy}{(x^2+y^2)^{\frac32}}, \frac{x^2}{(x^2+y^2)^{\frac32}}, 0\right).
\]
Errors and condition numbers for three mesh sizes, $h=\frac{1}{13}, \frac{1}{26}, \frac{1}{52},$ are shown in \Cref{fig:error_and_condition_plots}. 
A heat map of the $x$-component of the solution $\eta$ is shown in \Cref{fig:HL_ux}. We note that the expected rate of convergence $\mathcal O(h)$ is indeed achieved for the errors on the curl and gradient, with a slightly better rate for the fields themselves (going possibly towards $\mathcal O(h^2)$). The condition number for the scheme based on the ghost stabilisation remains at a reasonable magnitude, while for the unstabilised scheme (using the inner product $(\bullet,\bullet)_\Omega$ instead of $(\bullet,\bullet)_s$) it blows up to a level where solving the linear system becomes infeasible. 

\begin{figure}
  \centering
  \includegraphics[width=0.75\textwidth]{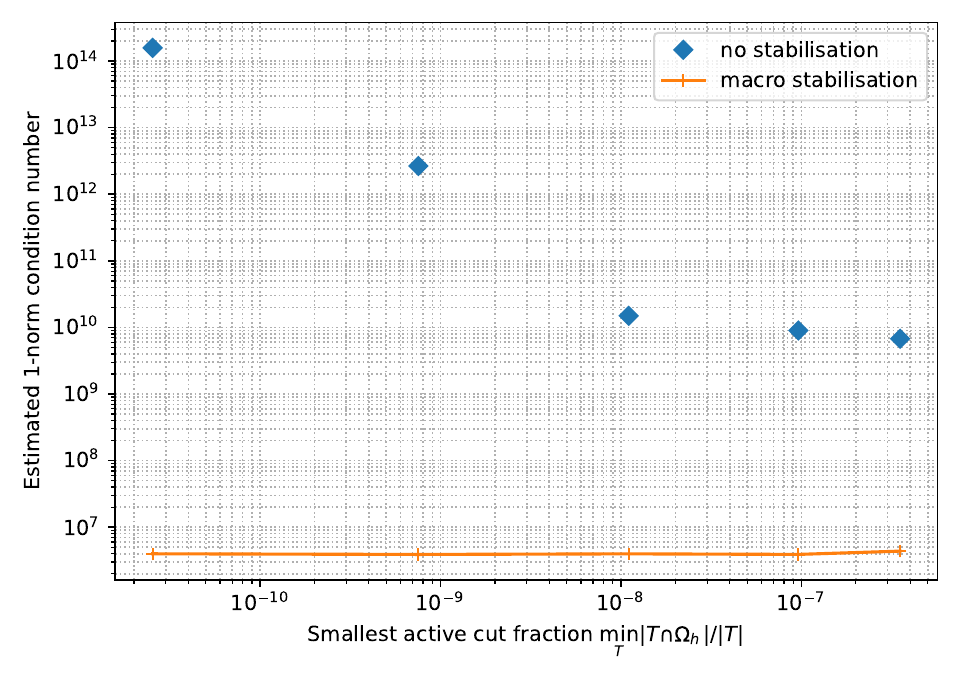}
  \caption{
    Hodge Laplace example: Condition numbers versus smallest cut fraction $\min_{T\in\mcT_h} |T\cap\Omega|/|T|$ for the stabilised and unstabilised system matrices. 
  }
  \label{fig:cut_position_condition}
\end{figure}

We end this example by illustrating with \Cref{fig:cut_position_condition} the independence of the condition number of the stabilised system matrix with respect to the relative position between the domain and the background mesh. The stabilised system matrix has a condition number that is independent of the cut position, while the unstabilised system matrix has a condition number that blows up as the cut fraction goes to zero.

\subsection{An $H^{\crl}$-conforming unfitted Stokes scheme}
Let all user parameters be set as in the previous example. 
Let $S_h^0 := P^-_1\Lambda^0\Omdh / \RR$ denote the continuous trimmed Lagrange space with zero mean.
With $f_h:=\Pi_h^{\crl}f$, the unfitted method for the Stokes equations in curl-curl formulation \cite{Girault:90,Bonelle:14, da2022presrobstokes}, with natural  boundary conditions $n \times \crl u_h = 0$ and $ n \cdot u_h = 0$ on $\partial\Omega$,
read: find
$(u_h,p_h)\in P^-_1\Lambda^1\Omdh\times S_h^0$ such that
\begin{equation}\label{eq:stokes.hcurl.unfitted}
    (\crl u_h,\crl v_h)_s
    +(\nabla p_h,v_h)_s
    +(u_h,\nabla q_h)_s
    =
    (f_h,v_h)_s\quad \forall (v_h,q_h)\in P^-_1\Lambda^1\Omdh\times S_h^0,
\end{equation}
where $f_h$ is the interpolation on $P^-_1\Lambda^1\Omdh$ of $f$. Note that this requires to define an extension of $f$ outside $\Omega$ that can be interpolated on this space, similar to \Cref{rem:poisson_div_constraint}. This definition of the load through interpolator and use of the stabilised product, as in the pressure-velocity coupling term, ensures the pressure-robustness of the scheme \cite[Remark 5]{da2022presrobstokes}.

\begin{figure}
  \centering
  \includegraphics[width=0.75\textwidth]{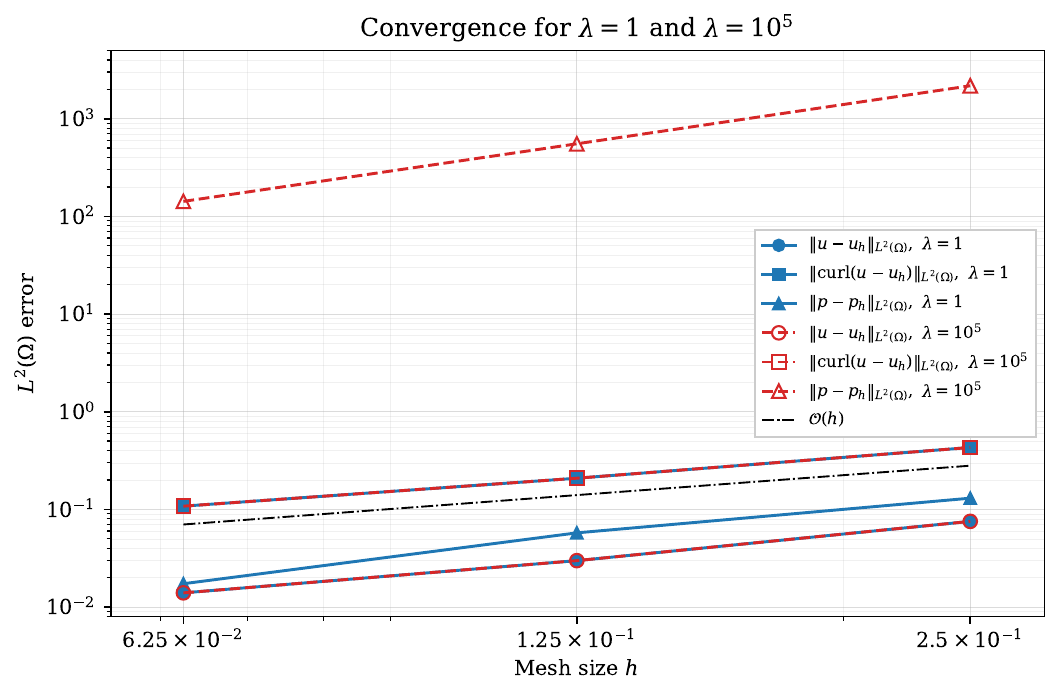}
  \caption{
    Stokes example: Convergence of the $L^2$-errors in the velocity, its curl, and the
    pressure for $\lambda=1$ and $\lambda=10^5$. Blue solid curves
    correspond to $\lambda=1$, while red dashed curves correspond to
    $\lambda=10^5$; the markers distinguish $u$, $\crl u$, and $p$.
    The velocity and curl errors coincide for the two pressure amplitudes,
    demonstrating pressure robustness of the unfitted scheme.
  }
  \label{fig:stokes.pressure.robustness}
\end{figure}

For $R=2/3$ we take $\Omega=B_R(0)\subset\mathbb [0,1]^3$ represented by a level set function, and set
\[
    s(x):=R^2-|x|^2,\qquad
    u(x):=6s(x)^2(-x_2,x_1,0),
    \qquad
    p_\lambda(x):=\lambda
    \left(|x|^2-\frac{3R^2}{5}\right).
\]
Then $\int_\Omega p_\lambda=0$, $\dive u=0$, and $u\cdot n = 0$ and $n\times (\crl u)=0$ on $\partial\Omega$. The right hand side is
\[
    f_\lambda
    :=
    \crl\crl u+\nabla p_\lambda
    =
    24\bigl(7|x|^2-5R^2\bigr)(x_2,-x_1,0)
    +2\lambda x.
\]
Theoretical estimates proved, e.g., in \cite{da2022presrobstokes}, show that the error on the velocity is independent of the pressure, and in particular of $\lambda$ in the chosen example. This is illustrated in \Cref{fig:stokes.pressure.robustness}, in which we can see that the errors on $u-u_h$ and $\crl u -\crl u_h$ remain exactly the same when $\lambda$ increases. Interestingly, if on the right hand side one puts $(f_h,v_h)_\Omega$, without stabilisation, the result is not pressure robust; see \Cref{tab:stokes.pressure.robustness}.

\begin{table}
  \centering
  \begin{tabular}{c | c c | c c}
    \hline
     & \multicolumn{2}{c|}{$\lambda=1$} & \multicolumn{2}{c}{$\lambda=10^5$} \\
    \hline
    $h$ & $\|u-u_h\|_{L^2}$ & $\|\mathrm{curl}\,u-\mathrm{curl}\,u_h\|_{L^2}$ & $\|u-u_h\|_{L^2}$ & $\|\mathrm{curl}\,u-\mathrm{curl}\,u_h\|_{L^2}$ \\
    \hline
    0.25   & 0.0756969 & 0.428067 & 226.69    & 2273.32 \\
    0.125  & 0.0299591 & 0.208658 & 30.0879   & 537.504 \\
    0.0625 & 0.013969  & 0.107864 & 4.59413   & 138.911 \\
    \hline
  \end{tabular}
  \caption{Stokes example: Non pressure-robust errors of the velocity and its curl for $\lambda=1$ and $\lambda=10^5$ when the load is $(f_h,v_h)_\Omega$ instead of $(f_h,v_h)_s$.}
  \label{tab:stokes.pressure.robustness}
\end{table}

\section*{Acknowledgements}

Funded by the European Union (ERC Synergy, NEMESIS, project number 101115663).
Views and opinions expressed are however those of the authors only and do not necessarily reflect those of the European Union or the European Research Council Executive Agency. Neither the European Union nor the granting authority can be held responsible for them.


\addcontentsline{toc}{section}{Bibliography}
\bibliographystyle{abbrvnat}
\bibliography{ref}

\end{document}